%% file: gquots.tex
\documentclass[a4paper]{amsart}
\usepackage{amssymb}
\usepackage[all]{xy}
\usepackage{amsmath}
\usepackage{graphics}
\usepackage{stmaryrd}
\usepackage{color}
\usepackage{dsfont}

\CompileMatrices
\newtheorem*{introthm}{Theorem}

\newtheorem{theorem}{Theorem}[section]
\newtheorem{lemma}[theorem]{Lemma}
\newtheorem{proposition}[theorem]{Proposition}
\newtheorem{corollary}[theorem]{Corollary}
\theoremstyle{definition}
\newtheorem{definition}[theorem]{Definition}
\newtheorem{example}[theorem]{Example}

\newtheorem{construction}[theorem]{Construction}

\newtheorem{remark}[theorem]{Remark}

\def\div{{\rm div}}

\def\rlv{{\rm rlv}}
\def\Eff{{\rm Eff}}
\def\Mov{{\rm Mov}}
\def\SAmple{{\rm SAmple}}
\def\Ample{{\rm Ample}}

\def\quot{/\!\!/}
\def\mal{\! \cdot \!}

\def\reg{{\rm reg}}
\def\rq#1{\widehat{#1}}
\def\t#1{\widetilde{#1}}
\def\b#1{\overline{#1}}
\def\bangle#1{\langle #1 \rangle}

\def\CC{{\mathbb C}}
\def\KK{{\mathbb K}}

\def\ZZ{{\mathbb Z}}

\def\QQ{{\mathbb Q}}
\def\PP{{\mathbb P}}
\def\XX{{\mathbb X}}

\def\WDiv{\operatorname{WDiv}}

\def\Cl{\operatorname{Cl}}

\def\Pic{\operatorname{Pic}}

\def\Hom{{\rm Hom}}
\def\GL{{\rm GL}}
\def\SO{{\rm SO}}

\def\Supp{{\rm Supp}}
\def\Spec{{\rm Spec}}
\def\Proj{{\rm Proj}}

\def\cone{{\rm cone}}
\def\lin{{\rm lin}}

\def\topto#1{\stackrel{{\scriptscriptstyle #1}}{\longrightarrow}}

\def\cov{{\rm cov}}
\def\SL{{\rm SL}}
\def\Sp{{\rm Sp}}

\newcounter{itemnumber}

\begin{document}

\sloppy

\title[GIT via Cox rings]
{Geometric Invariant Theory \\ via Cox rings}

\author[I.~Arzhantsev]{Ivan~V. Arzhantsev} 
\thanks{Supported by INTAS YS 05-109-4958}
\address{Department of Higher Algebra, 
Faculty of Mechanics and Mathematics, 
Moscow State Lomonosov University,
Leninskie Gory, GSP-2, Moscow, 119992, Russia}
\email{arjantse@mccme.ru}
\author[J.~Hausen]{J\"urgen Hausen} 
\address{Mathematisches Institut, Universit\"at T\"ubingen,
Auf der Morgenstelle 10, 72076 T\"ubingen, Germany}
\email{hausen@mail.mathematik.uni-tuebingen.de}
\subjclass[2000]{14L24, 14L30, 14C20}

\begin{abstract}
We consider actions of reductive groups
on a variety with finitely generated 
Cox ring, e.g., the classical case of 
a diagonal action on a product of projective spaces.
Given such an action, we construct via combinatorial data
in the Cox ring all maximal open subsets such 
that the quotient is quasiprojective or embeddable 
into a toric variety.
As applications, we obtain an explicit description
of the chamber structure of the linearized ample 
cone and several Gelfand-MacPherson type 
correspondences relating quotients by reductive
groups to quotients by torus actions. 
Moreover, our approach provides a general access
to the geometry of many of the resulting quotient
spaces. 
\end{abstract}

\maketitle

\section{Introduction}

The passage to a quotient by an algebraic group 
action is often an essential step in classical moduli
space constructions of Algebraic Geometry,
and it is the task of Geometric Invariant Theory
(GIT) to provide such quotients.
Starting with Mumford's approach of constructing
quotients for actions of reductive groups
on projective varieties via
linearized  line bundles and their sets of
semistable points~\cite{Mu}, the notion
of a ``good quotient'' became a central concept
in GIT, compare~\cite{Se} and~\cite{BB2}.
Recall that a good quotient for an action 
of a reductive group $G$ on a variety
$X$ is an affine morphism 
$\pi \colon X \to Y$ of varieties
such that $Y$ carries the sheaf of invariants 
$(\pi_* \mathcal{O}_X)^G$
as its structure sheaf.
In general, a $G$-variety $X$ need not admit
a good quotient, but there may be many different
invariant open subsets $U \subseteq X$ with a good
quotient; we will call them the good $G$-sets.

In this paper, we consider $G$-varieties $X$ with 
a finitely generated Cox ring, e.g. $X$ being a 
product of projective spaces, and ask for 
good $G$-sets $U \subseteq X$, which are maximal 
with respect to the properties either that the 
quotient space $U \quot G$ is quasiprojective or, 
more generally, that it comes with the A2-property;
the latter means that any two points of $U \quot G$ 
admit a common affine neighbourhood, or, equivalently,
that $U \quot G$ admits a closed embedding into 
some toric variety, see~\cite{Wl}.
Our aim is to provide a constructive approach to such
good $G$-sets, splitting the explicit computation 
into two parts: firstly computations of invariant 
rings in the spirit of classical Invariant Theory 
and, secondly, combinatorial computations with
convex polyhedral cones.
Another feature is that our approach opens an access
to the geometry of quotient spaces via the methods
developed in~\cite{BeHa1}.

Let us present our results more in detail.
A first step is to consider actions of $G$ on 
factorial affine varieties $X$.
The basic data for the construction
of good $G$-sets of $X$ 
are {\em orbit cones\/}. They live 
in the rational character space $\XX_\QQ(G)$,
and for any $x \in X$ its orbit cone $\omega(x)$ is the 
convex cone generated by all $\chi \in \XX(G)$
admitting a semiinvariant $f$ with weight $\chi$ 
such that $f(x) \ne 0$ holds.
It turns out that there are only finitely many 
orbit cones and all of them are polyhedral.

\goodbreak

Based on the concept of orbit cones, we introduce 
the data describing the good $G$-sets of the 
factorial affine variety $X$.
First, we associate to any character 
$\chi \in \XX(G)$ its {\em GIT-cone}, namely
$$
\lambda(\chi) 
\ := \ 
\bigcap_{\chi \in \omega(x)}  \omega(x)
\ \subseteq \
\XX_\QQ(G).
$$
Secondly, we say that a collection $\Phi$ 
of orbit cones is {\em 2-maximal}, 
if for any two members their 
relative interiors overlap and $\Phi$   
is maximal with respect to this property.
Here comes the first result,
see Theorems~\ref{thm:GITfan} and~\ref{thm:maxcoll}.

\begin{introthm}
Let a connected reductive group $G$ act on a 
factorial affine variety~$X$.
\begin{enumerate}
\item
The GIT-cones form a fan in $\XX_\QQ(G)$,
and this fan is in a canonical order reversing bijection
with the collection of sets of semistable 
points of~$X$.
\item
There is a canonical bijection from the set of 2-maximal 
collections of orbit cones onto the collection of 
A2-maximal good $G$-sets of $X$.
\end{enumerate}
\end{introthm}

For the case of a torus $G$ this result was known before.
The first statement is given in~\cite{BeHa2}. 
Moreover, a result similar to the second statement 
was obtained in~\cite{BBSw5}
for linear torus actions on vector spaces, 
and for torus actions on 
any affine factorial $X$, 
statement~(ii) is given in~\cite{ArHa}.

To obtain the general statement, we 
reduce to the case of a torus action as follows.
Consider the quotient $Y := X \quot G^s$
by the semisimple part $G^s \subseteq G$.
It comes with an induced action of the 
torus $T := G/G^s$, and the key observation is 
that the good $T$-sets in $Y$
are in a canonical bijection with the  
good $G$-sets in $X$, see Proposition~\ref{thm:reduct}.
Note that this is the place, where in explicit computations, 
Classical Invariant Theory comes in, as it provides 
often the necessary information on the algebra $\KK[X]^{G^s}$ 
of invariants, see the examples discussed
in Sections~\ref{sec:amplegit} and~\ref{sec:quotgeo}.

The second step is passing to the case 
of a normal variety $X$ with a finitely 
generated Cox ring $\mathcal{R}(X)$;
recall that, for its definition,
one assumes that the divisor class 
group $\Cl(X)$ is free and finitely generated,
and then sets
\begin{eqnarray*}
\mathcal{R}(X)
& := & 
\bigoplus_{D \in \Cl(X)} \Gamma(X,\mathcal{O}(D)).
\end{eqnarray*}
The ``total coordinate space'' $\b{X}$ of $X$
is the spectrum of the Cox ring 
$\mathcal{R}(X)$.
This  $\b{X}$ is a factorial 
affine variety, see~\cite{BeHa2},
acted on by the Neron-Severi torus $H$
having the divisor class group $\Cl(X)$
as its character lattice.
Moreover, $X$ can be reconstructed from 
$\b{X}$ as a good quotient 
$q \colon \rq{X} \to X$ by $H$ for an
open subset $\rq{X} \subseteq \b{X}$,
see Section~\ref{sec:lift} for details.

After replacing $G$ with a simply connected 
covering group, its action on $X$ 
can be lifted to
the total coordinate space $\b{X}$.
The actions of $H$ and $G$  on $\b{X}$
commute, and thus 
define an action of the direct product
$\b{G} := H \times G$.
Given a good $\b{G}$-set $W \subseteq \b{X}$,
we introduce in~\ref{def:satint} a 
``saturated intersection''
$W \sqcap_G \rq{X}$.
The main feature of this construction is
the following, see Theorem~\ref{thm:sqcap}.

\begin{introthm}
The canonical assignment
$W \mapsto q(W \sqcap_G \rq{X})$ 
defines a surjection from the 
collection of good $\b{G}$-sets in $\b{X}$ 
to the collection of good $G$-sets in $X$. 
\end{introthm}

So this result reduces the construction of 
good $G$-sets on $X$ to the construction 
of good $\b{G}$-sets in $\b{X}$,
and the latter problem, as noted before,
is reduced to the case of a torus action.
Again, this allows explicit computations.
Note that our way to reduce the construction
of quotients to the case of a torus action
has nothing in common with the various 
approaches based on the Hilbert-Mumford 
Criterion, see~\cite{BB2}, \cite{DoHu}, 
\cite{Mu}, \cite{Re}, but is rather in the 
spirit of~\cite[Sec.~3]{Tha}.

\goodbreak

As a first application of this result, 
we give an explicit description of  
the ample GIT-fan, i.e., the chamber structure 
of the linearized ample cone,
for a given normal projective 
$G$-variety $X$ with finitely generated 
Cox ring, see Proposition~\ref{amplegit}; 
recall that existence
of the ample GIT-fan for any normal projective 
$G$-variety was proven in~\cite{DoHu}~\cite{Tha},
and, finally,\cite{Re}. 
As an example, we compute the ample GIT-fan for 
the diagonal action of $\Sp(2n)$ on a product of 
projective spaces $\PP^{2n-1}$, see~Theorem~\ref{tsp}.

As a second application of the above result 
we obtain Gelfand-MacPherson type correspondences.
Classically~\cite{GM}, this correspondence 
relates orbits of the diagonal action of the special 
linear group $G$ on a product of projective spaces
to the orbits of an action of a torus $T$ on a 
Grassmannian. 
Kapranov~\cite{Ka} extended this  
correspondence to isomorphisms of certain 
GIT-quotients and used it in his study 
of the moduli space of point configurations 
on the projective line.
Similarly, Thaddeus~\cite{Tha2} proceeded
with complete collineations.
In Section~\ref{sec:gmcorr}, we put these correspondences 
into a general framework, relating GIT-quotients 
and also their inverse limits. 
As examples, we retrieve a result of~\cite{Tha2} 
and also an isomorphism of GIT-limits in the 
setting of~\cite{Ka}.

Finally, we use our approach to study the geometry 
of quotient spaces of a connected reductive 
group $G$ on a normal variety $X$ 
with finitely generated Cox ring. 
The basic observation is that in many cases
our quotient construction provides  
the Cox ring of the quotient spaces.
This allows to apply the language of 
bunched rings developed in~\cite{BeHa2},
which encodes information on the geometry
of a variety in terms of combinatorial data
living in the divisor class group.

\tableofcontents

\section{Some background on good quotients}
\label{sec:goodquot}

In this section, we recall the concept of a good 
quotient and state basic properties, which 
will be used freely in the subsequent text.

Throughout the whole paper,
we work in the category of algebraic varieties 
over an algebraically closed field $\KK$ of
characteristic zero.
By a point we always mean a closed point.
If we say that an algebraic group $G$ acts on a 
variety $X$, then we tacitly assume that this 
action is given by a morphism $G \times X \to X$,
and we refer to $X$ as a $G$-variety.
As usual, we say that a morphism 
$\varphi \colon X \to Y$ of $G$-varieties is 
equivariant if it is compatible with the 
actions in the sense that always
$\varphi(g \mal x) = g \mal \varphi(x)$ holds.
Moreover, a morphism is called invariant, if
it is constant along the orbits.

The classical finiteness theorem in Invariant Theory
says that for an
action of a reductive linear algebraic group 
on an affine 
variety $X = \Spec(A)$, the algebra $A^G$ 
of invariant functions is finitely generated. 
This allows to define the classical invariant
theory quotient $Y := \Spec(A^G)$, which comes 
with a morphism $p \colon X \to Y$.
The notion of a good quotient is locally modeled
on this concept:

\begin{definition}
Let $G$ be a reductive linear algebraic group.
A {\em good quotient\/} for a 
$G$-variety $X$ is an affine 
morphism $p \colon X \to Y$
onto a variety $Y$ 
such that the pullback 
$p^* \colon \mathcal{O}_Y \to (p_* \mathcal{O}_X)^G$
to the sheaf of invariants 
is an isomorphism.
A good quotient is called {\em geometric}, 
if its fibers are precisely
the orbits. 
\end{definition}


The basic properties of a good quotient 
$p \colon X \to Y$ of a $G$-variety are
that it sends closed $G$-invariant 
subsets $A \subseteq X$ to closed sets
$p(A) \subseteq Y$, and that for any
two disjoint closed $G$-invariant subsets
$A,A' \subseteq X$ their images 
$p(A), p(A') \subseteq Y$ 
are again disjoint.
An immediate consequence is 
that each fiber $p^{-1}(y)$
of a good quotient $p \colon X \to Y$ 
contains precisely one closed $G$-orbit,
and this orbit lies in the closure of any further 
orbit in $p^{-1}(y)$.

These basic properties imply 
that a good quotient $X \to Y$ 
for a $G$-variety $X$ is categorical, 
i.e., any $G$-invariant morphism $X \to Z$ 
factors uniquely through $X \to Y$. 
In particular, good quotient spaces are 
unique up to 
isomorphism. 
This justifies the notation
$X \to X \quot G$ for good and $X \to X/G$ 
for geometric quotients, which we will use 
frequently in the sequel.

\begin{proposition}
\label{goodquotprops}
Let $G$ be a  connected reductive group,
$H \subseteq G$ a normal, reductive subgroup,
and $X$ be a $G$-variety.
\begin{enumerate}
\item
If the good quotient $X \to X \quot H$ exists, 
then there is a unique $G$-action on $X \quot H$ 
making  $X \to X \quot H$ equivariant, and 
this action uniquely induces an action
of $G/H$ on $X \quot H$.
\item
The good quotient $X \to X \quot G$ exists if 
and only if the good quotients $X \to X \quot H$
and $X \quot H \to (X \quot H) \quot (G/H)$ exist.
In this case, one has a commmutative diagram
$$ 
\xymatrix{
X 
\ar[rr]^{\quot H}
\ar[d]_{\quot G}
& &
X \quot H
\ar[d]^{\quot (G/H)}
\\
X \quot G \ar@{<->}[rr]_{\cong}
& &
(X \quot H) \quot (G/H)
}
$$
\end{enumerate}
\end{proposition}

\begin{proof}
In the setting of~(i),
universality of the good quotient
allows to push down the $G$-action to $X \quot H$,
see~\cite[Thm.~7.1.4]{BB2}.
In the setting of~(ii),
if $X \to X \quot G$ exists,
then also $X \to X \quot H$ 
exists, see~\cite[Cor.~10]{BBSw4},
and one directly verifies that
the induced morphism 
$X \quot H \to X \quot G$ 
is a good quotient for the action of
$G/H$.
Conversely, if the stepwise good 
quotients exist, then one directly 
verifies that their composition
is a good quotient for the $G$-variety
$X$.
\end{proof}

In general, quite a few non-affine $G$-varieties $X$
admit a good quotient $X \to X \quot G$.
However, there may be many open invariant subsets
$U \subset X$ with a good quotient $U \to U \quot G$,
and it is one of the main tasks of the theory of 
good quotients to describe all these sets.
Here we fix the basic terminology.

\begin{definition}
Let a connected reductive group $G$ 
act on a variety~$X$.
\begin{enumerate}
\item
By a {\em good $G$-set\/} we mean an open, $G$-invariant 
subset $U \subseteq X$ admitting a good quotient
$U \to U \quot G$.
\item
If $U \subseteq X$ is a good $G$-set,
then the {\em $G$-limit of $x \in U$\/}
is the unique closed $G$-orbit in the 
closure of $G \mal x$ with respect to $U$;
we denote it by $\lim_{G}(x,U)$.
\item
For a good $G$-set $U \subseteq X$,
we say that $U' \subseteq U$ is a 
{\em $G$-saturated inclusion\/} if 
$U'$ is open, $G$-invariant, and 
for any $x \in U'$ one has 
$\lim_G(x,U) \subseteq U'$. 
\end{enumerate}
\end{definition}

Given a good $G$-set $U \subseteq X$,
we have mutually inverse bijections
between the collection of 
$G$-saturated subsets $U' \subseteq U$ 
and the collection of open subsets 
$V \subseteq U \quot G$,
sending  $U' \subseteq U$ to 
$p(U') \subseteq U \quot G$
and 
$V \subseteq U \quot G$ to $p^{-1}(V)$.
Moreover, for any $G$-saturated 
$U' \subseteq U$, the restriction 
$p \colon U' \to p(U')$ is 
a good quotient for the $G$-variety $U'$.

The preceding observation allows 
to concentrate in the study of good $G$-sets
to certain maximal ones, where one 
also may impose properties on the quotient 
spaces, like quasiprojectivity or 
the A2-property, i.e.,  
any two points admit a common affine 
neighbourhood.
Here are the precise notions,
compare~\cite{BB2}.

\begin{definition}
Let a connected reductive group $G$ 
act on a variety~$X$.
We say that a good $G$-set $U \subseteq X$ 
is 
{\em maximal (qp-maximal, A2-maximal)\/},
if it is maximal with respect to
saturated inclusion among all good 
$G$-sets $W \subseteq X$ (among those
with $W \quot G$ quasiprojective, 
having the  A2-property).
\end{definition}

We conclude this section with recalling the 
construction of good $G$-sets with quasiprojective
quotient spaces as sets of semistable points
presented by D.~Mumford in~\cite{Mu}.
In fact, we use here a slightly more general 
version, based on Weil divisors instead of 
line bundles, see~\cite{Ha}; 
compared to Mumford's original approach 
this has the advantage to produce 
all qp-maximal subsets, see Prop.~\ref{qpquot}.

Let $X$ be a normal $G$-variety, 
where $G$ is a reductive linear
algebraic group. 
To any Weil divisor $D$ on $X$, 
we associate a sheaf of 
$\mathcal{O}_X$-algebras, 
and consider the corresponding
relative spectrum with its 
canonical morphism, compare~\cite[I.9.4.9]{EGA} 
and \cite[pp. 128/129]{Har}:
$$
\mathcal{A}
\ := \
\bigoplus_{n \in \ZZ_{\ge 0}}
\mathcal{O}_X(nD),
\qquad
X(D) \ := \
\Spec_X(\mathcal{A}),
\qquad
q_D \colon X(D) \to X.
$$
The $\ZZ_{\ge 0}$-grading of the 
sheaf of algebras $\mathcal{A}$ 
defines a $\KK^*$-action on  
$X(D)$ having the canonical morphism 
$q_D \colon X(D) \to X$ as a good 
quotient.
For these constructions, we tacitly 
assumed that $\mathcal{A}$ is locally 
of finite type over $\mathcal{O}_X$, and thus
$X(D)$ is in fact a variety; this
holds for all Cartier divisors $D$,
and will later be guaranteed by 
a global finiteness condition on $X$.

\begin{definition}
A {\em $G$-linearization\/} of the divisor 
$D$ is a (morphical) $G$-action on $X(D)$ 
that commutes with the $\KK^*$-action on $X(D)$
and makes $q_D \colon X(D) \to X$ into a 
$G$-equivariant morphism.
\end{definition}

Note that for a Cartier divisor $D$, the sheaf
$\mathcal{A}$ is locally free of rank one, 
and hence $X(D) \to X$ is a line bundle.
So, in this context a $G$-linearization 
of $D$ is a fiberwise linear action on the
total space $X(D)$ making the projection
$X(D) \to X$ equivariant; this is precisely 
Mumford's original notion of a linearization
of a line bundle.

On the (invariant) set $X_{\reg} \subseteq X$ 
of smooth points, 
any $G$-linearized Weil divisor $D$ is Cartier,
and hence $X(D) \to X$ is a $G$-linearized 
line bundle over $X_{\reg}$ in the usual sense.
This allows to define the (linearized) sum 
$D + D'$ of  linearized divisors $D,D'$ 
by extending the canonical action 
on $X(D) \otimes X(D')$ from 
$X_{\reg}$ to all of $X$, 
compare~\cite[Sec.~1]{BeHa2}.

There is also the concept of the linearized
divisor class group $\Cl_G(X)$.
Call two linearized divisors 
$D, D'$ equivalent, if there is an 
$(\KK^* \times G)$-equivariant isomorphism
$X(D) \to X(D')$. Then the addition defined just 
before induces a well defined addition on the
set  $\Cl_G(X)$ of classes of linearized 
divisors turning it into a group.
Note that we have a canonical restriction 
isomorphism $\Cl_G(X) \to \Pic_G(X_{\reg})$
to the linearized Picard group of $X_{\reg}$.

Finally, we come to the definition of semistability.
Given a $G$-linearized Weil divisor $D$,
one has a natural rational 
$G$-representation on the vector space of 
its global sections
$\Gamma(X,\mathcal{O}(D)) = \Gamma(X(D),\mathcal{O})$, 
namely
$$
G \times \Gamma(X(D), \mathcal{O}) 
\to 
\Gamma(X(D),\mathcal{O}),
\qquad
(g \mal f)(x) := f(g^{-1} \mal x).
$$ 
In particular, this allows to speak
about the space $\Gamma(X(D),\mathcal{O})^G$ 
of invariant sections of $D$.
Moreover, for any section 
$f \in \Gamma(X,\mathcal{O}(D))$,
one defines its set of zeroes $Z(f) \subseteq X$ 
by setting 
$$
Z(f)
\ := \
\Supp(\div(f)+D). 
$$


\begin{definition}
\label{semistabdef}
Let $G$ be a linear algebraic group,
$X$ a normal $G$-variety and 
$D$ a $G$-linearized Weil divisor
on $X$. 
\begin{enumerate}
\item
We call $x \in X$ 
{\em semistable\/} with respect to $D$
if there are 
$n \in \ZZ_{>0}$ and 
$f \in \Gamma(X,\mathcal{O}(nD))^G$ 
such that $X \setminus Z(f)$ is an
affine neighbourhood of $x$.
\item
The set of semistable points
of a $G$-linearized Weil divisor $D$ 
on $X$ will be denoted by $X^{ss}(D)$, 
or $X^{ss}(D,G)$, if the group $G$ 
needs to be specified.
\item
If $D'$ is another $G$-linearized 
Weil divisor on $X$, then we say that 
$D$ and $D'$ are {\em GIT-equivalent\/} 
if their associated sets of semistable 
points coincide.
\end{enumerate}
\end{definition}

Note that two linearized divisors defining 
the same class in $\Cl_G(X)$ have the same 
set of semistable points.
From~\cite[Theorem~3.3]{Ha} we infer the 
following features of the 
sets of semistable points:

\begin{proposition}
\label{qpquot}
Let $G$ be a reductive linear algebraic group
and $X$ a normal $G$-variety.
\begin{enumerate}
\item If $D$ is a $G$-linearized Weil divisor on $X$,
then there exists a good quotient 
$X^{ss}(D) \to X^{ss}(D) \quot G$ with a 
quasiprojective quotient space. 
\item If $U \subset X$ is a $G$-invariant
open subset having a good quotient 
$U \to U \quot G$ with $U  \quot G$ 
quasiprojective, then $U$ is $G$-saturated
in some set $X^{ss}(D)$.
\end{enumerate}
\end{proposition}


\section{Good quotients of factorial affine varieties}
\label{sec:afffact}

In this section, we consider an action of a
connected reductive group $G$ on an irreducible
affine variety $Z$.
In the first result, we describe the 
collection of sets of semistable points 
arising from the possible linearizations 
of the trivial line bundle over $Z$,
and in the second one, we describe the collection
of A2-maximal subsets of $Z$ provided 
that $Z$ is factorial.
Both descriptions are of combinatorial nature and are 
given in terms of
certain convex polyhedral cones.
The first setting was also studied 
in~\cite{Ki}; there a numerical criterion for
semistability was given.

Let us briefly fix the necassary notation.
Given a polyhedral cone $\sigma$ 
in some rational vector space, we denote by 
$\sigma^\circ \subseteq \sigma$ its relative 
interior, and for $\tau \subseteq \sigma$, 
we write 
$\tau \preceq \sigma$ if $\tau$ is a face 
of $\sigma$.
By a fan in a rational vector space,
we mean a finite collection
$\Sigma$ of convex, polyhedral cones such that 
for $\sigma \in \Sigma$ also every 
face $\tau \preceq \sigma$ 
belongs to $\Sigma$ and for any two 
$\sigma_1, \sigma_2 \in \Sigma$ one has 
$\sigma_1 \cap \sigma_2 \preceq \sigma_i$;
note that we don't require the cones of $\Sigma$ 
to be pointed.

Now we turn to the $G$-linearizations 
of the trivial bundle $Z \times \KK \to Z$;
they arise from the elements $\chi \in M$ 
of the character group $M := \XX(G)$ as 
follows.
\begin{eqnarray} 
\label{eqn:lin2char}
G \times (Z \times \KK) \ \to \ Z \times \KK,
\qquad
g \mal (z,z') \ = \ (g \mal z, \chi(g) z').
\end{eqnarray}
Every such $G$-linearization defines a set 
$Z^{ss}(\chi) \subseteq Z$ of semistable points, and this 
set is explicitly given by 
\begin{eqnarray*}
Z^{ss}(\chi)
& = & 
\{z \in Z; \; f(z) \ne 0 
\text{ for some } 
f \in \Gamma(Z,\mathcal{O})_{n\chi}, \; n > 0 \}.
\end{eqnarray*}
As outlined before, the set $Z^{ss}(\chi)$ admits a 
good quotient for the action of $G$; the quotient
space is given by 
$$ 
Z^{ss}(\chi) \quot G \ = \ \Proj(A(\chi)),
\quad \text{where }
A(\chi) \ := \ \bigoplus_{n \in \ZZ_{\ge 0}}  \Gamma(Z,\mathcal{O})_{n\chi}.
$$
In particular, 
$Z^{ss}(\chi) \quot G$ is projective over $Z \quot G$.
Our description
of the collection of sets $Z^{ss}(\chi)  \subseteq Z$ 
is formulated in terms of the following 
combinatorial data.


\begin{definition}
Let $G$ be a  connected reductive group, denote by
$M_\QQ := M \otimes_\ZZ \QQ$ its rational character 
space, and let 
$Z$ be an irreducible affine $G$-variety.
\begin{enumerate}
\item
The {\em weight cone\/} of the $G$-variety $Z$ is the 
convex cone $\omega(Z) \subseteq M_\QQ$ generated by 
all $\chi \in M$ with $\Gamma(Z,\mathcal{O})_\chi \ne 0$.
\item
The {\em orbit cone\/} of a point $z \in Z$ is the 
convex cone $\omega(z) \subseteq M_\QQ$ generated by 
all $\chi \in M$ that admit an 
$f \in \Gamma(Z,\mathcal{O})_\chi$ with $f(z) \ne 0$.
\item
The {\em GIT-cone\/} $\lambda(\chi) \subseteq M_\QQ$ 
of a character $\chi \in M$ is the intersection of 
all orbit cones containing $\chi$:
\begin{eqnarray*}
\lambda(\chi) 
& := & 
\bigcap_{\chi \in \omega(x)}  \omega(x).
\end{eqnarray*}
\item 
The {\em GIT-fan\/} of the $G$-variety $Z$ is the 
collection $\Sigma(Z)$ of all GIT-cones 
$\lambda(\chi)$, where $\chi \in M$.
\end{enumerate}
\end{definition}

\begin{theorem}
\label{thm:GITfan}
Let $G$ be a  connected reductive group and 
$Z$ be an irreducible affine $G$-variety.
\begin{enumerate}
\item
The weight cone, the orbit cones and the GIT-cones 
of the $G$-action on $Z$ are all polyhedral, 
and there are only finitely many of them.
\item
For every character $\chi \in M$, 
the associated set of semistable points $Z^{ss}(\chi)  \subseteq Z$ 
is given by
\begin{eqnarray*}
Z^{ss}(\chi)
& = & 
\{z \in Z; \; \chi \in \omega(z) \}.
\end{eqnarray*}
\item
The GIT-fan $\Sigma(Z)$ is in fact a fan 
in the rational character space 
$M_\QQ$, and the union of all 
$\lambda(\chi) \in \Sigma(Z)$  
is precisely the weight cone $\omega(Z)$.
\item
For any $\chi \in M$, the set $Z^{ss}(\chi)  \subseteq Z$  
is nonempty if and only if $\chi \in \omega(Z)$. 
Moreover, for any two $\chi,\chi' \in \omega(Z) \cap M$,
one has
\begin{eqnarray*}
Z^{ss}(\chi) \ \subseteq \ Z^{ss}(\chi')
& \iff & 
\lambda(\chi) \ \succeq \ \lambda(\chi').
\end{eqnarray*}
\item
If $Z$ is factorial, then  $\Sigma(Z)$ is in 
bijection to the qp-maximal good $G$-sets of $Z$ via
$\lambda \mapsto Z^{ss}(\chi)$, with any $\chi$
taken from the relative interior of $\lambda$.
\end{enumerate}
\end{theorem}

We prove these assertions by reducing them to the
known case of a torus action.
For this, we consider the semisimple part 
of $G$, i.e., 
the maximal connected semisimple subgroup 
$G^s \subseteq G$.
Recall that $G^s \subseteq G$ is a 
normal subgroup,
the factor group $T := G/G^s$ is a torus,
and $G \to T$ induces an isomorphism of 
the character groups.
The latter allows us to identify 
the character groups of $G$ and $T$;
we denote them both by $M$.

\begin{lemma}
\label{lem:reduct}
Let $G$ be a  connected reductive group and
$Z$ an irreducible affine 
$G$-variety.
Consider the quotient 
$\pi \colon Z \to Z \quot G^s$
and the induced action 
$T := G/G^s$ on $Y := Z \quot G^s$.
\begin{enumerate}
\item
For every point $z \in Z$, we have $\omega(\pi(z)) = \omega(z)$.
\item
For every character $\chi \in M$, we have 
$Z^{ss}(\chi) = \pi^{-1}(Y^{ss}(\chi))$.
\item
If $Z$ is factorial, then also
$Z \quot G^s$ is factorial.
\end{enumerate}
\end{lemma}

\begin{proof}
The first two assertions follow directly from 
the fact that 
the algebra $\Gamma(Z, \mathcal{O})^{G^s}$ of 
$G^s$-invariants in  $\Gamma(Z, \mathcal{O})$
equals the algebra of $G$-semiinvariants 
of $\Gamma(Z, \mathcal{O})$. 
The third one is well known, see~\cite[Thm.~3.17]{PV2}.
\end{proof}

\begin{proof}[Proof of Theorem~\ref{thm:GITfan}]
We first turn to statements~(i) to~(iv).
By Lemma~\ref{lem:reduct} it suffices to 
have the corresponding statements for the action
of the torus $T := G/G^s$ on the affine variety
$Y := Z \quot G^s$.
In that case, the statements were proven 
in~\cite[2.5, 2.7, 2.9 and~2.11]{BeHa2}.

To see~(v), recall first from Prop.~\ref{qpquot}
that every good $G$-set $W \subseteq Z$ with 
$W \quot G$  quasiprojective
is $G$-saturated in the set of semistable 
points $Z^{ss}(D)$ of some linearized divisor.
Since $Z$ is factorial, we have 
$Z^{ss}(D) = Z^{ss}(\chi)$ for some 
$\chi \in M$.
This consideration shows in particular that
every qp-maximal good $G$-set
$W \subseteq Z$ is of the form $W = Z^{ss}(\chi)$ 
for some $\chi \in M$.

So, in order to establish the bijection as claimed,
we only have to show that any 
$Z^{ss}(\chi) \subseteq Z$ is qp-maximal.
Suppose that some $Z^{ss}(\chi) \subseteq Z$ is not. 
Then we have a $G$-saturated 
inclusion $Z^{ss}(\chi) \subseteq Z^{ss}(\chi')$ 
with a qp-maximal $Z^{ss}(\chi') \subseteq Z$
and a commutative diagram
$$ 
\xymatrix{
Z^{ss}(\chi) 
\ar[rr] \ar[d]
& &
Z^{ss}(\chi')
\ar[d]
\\
Z^{ss}(\chi) \quot G
\ar[rr] \ar[dr]
& &
Z^{ss}(\chi') \quot G
\ar[dl]
\\
& Z \quot G &
}
$$
Since the quotient spaces in the 
middle line are projective
over $Z \quot G$, the induced 
morphism 
$Z^{ss}(\chi) \quot G \to 
Z^{ss}(\chi') \quot G$
is projective. On the other hand,
by $G$-saturatedness, it is an open
embedding.
This implies
$Z^{ss}(\chi) \quot G = 
Z^{ss}(\chi') \quot G$,
and thus, again by 
$G$-saturatedness,
$Z^{ss}(\chi) = Z^{ss}(\chi')$.
\end{proof}

Our next aim is a description of all A2-maximal
good $G$-sets of a factorial affine $G$-variety 
$Z$. 
The necessary data are again given in terms of 
orbit cones.

\begin{definition}
Let $G$ be a connected  reductive group and  
$Z$ an irreducible affine $G$-variety.
Let $\Omega(Z)$ denote the collection of
all orbit cones $\omega(z)$, where $z \in Z$.
\begin{enumerate}
\item
By a {\em 2-connected collection\/}
we mean a subcollection
$\Psi \subseteq \Omega(Z)$ 
such that
$\tau_1^\circ \cap \tau_2^\circ
\ne \emptyset$ 
holds for any two $\tau_1, \tau_2 \in \Psi$.
\item
By a {\em 2-maximal collection\/}, 
we mean a 2-connected collection, which is not 
a proper subcollection of any other 
2-connected collection.
\item
We say that a 2-connected collection $\Psi$ is a 
{\em face\/} of a 2-connected collection $\Psi'$ 
(written $\Psi \preceq \Psi')$, if for any 
$\omega' \in \Psi'$ there is an $\omega \in \Psi$
with $\omega \preceq \omega'$. 
\item
To any collection $\Psi \subseteq \Omega(Z)$,
we associate a $G$-invariant subset
$U(\Psi) \subseteq Z$ as follows:
\begin{eqnarray*}
U(\Psi)
& := & 
\{ z \in Z; \; \omega_0 \preceq \omega(z) 
\text{ for some } \omega_0 \in \Psi \}.
\end{eqnarray*} 
\item
To any $G$-invariant open subset $U \subseteq Z$, 
we associate a set of orbit cones, namely
\begin{eqnarray*}
\Psi(U)
& := &  
\{\omega(z); \; z \in U \text{ with } 
G \mal z \text{ closed in } U
\}.
\end{eqnarray*}
\end{enumerate}
\end{definition}

\begin{theorem}
\label{thm:maxcoll}
Let $G$ be a  connected reductive group, and let 
$Z$ be a factorial affine $G$-variety.
Then we have mutually inverse 
bijections of finite sets:
\begin{eqnarray*}
\left\{
 \text{2-maximal collections in } \Omega(Z)
\right\}
& \longleftrightarrow &
\{\text{A2-maximal good $G$-sets of } Z \}
\\
\Psi 
& \mapsto & 
U(\Psi)
\\
\Psi(U)
& \mapsfrom &
U
\end{eqnarray*}
These bijections are order-reversing maps of 
partially ordered sets in the sense 
that we always have 
\begin{eqnarray*}
\Psi \ \preceq \ \Psi'
& \iff &
U(\Psi) \ \supseteq \ U(\Psi').
\end{eqnarray*}
\end{theorem}

As before, the idea of proof is to reduce the
problem via passing to the quotient $Z \quot G^s$ 
to the case of a torus action.
This time the reduction step is a statement 
of independent interest.

\begin{proposition}
\label{thm:reduct}
Let $G$ be a connected reductive group and
$Z$ a factorial affine $G$-varie\-ty.
Consider the semisimple part 
$G^s \subseteq G$, the torus $T := G/G^s$,
the quotient $\pi \colon Z \to Z \quot G^s$,
and the induced $T$-action on $Z \quot G^s$.
Then we have mutually inverse bijections
\begin{eqnarray*}
\left\{
\text{good $G$-sets of $Z$}
\right\}
& \longleftrightarrow &
\left\{
\text{good $T$-sets of $Z \quot G^s$}
\right\}
\\
U & \mapsto & \pi(U), \\
\pi^{-1}(V) & \mapsfrom  & V.
\end{eqnarray*}
Both assignments preserve saturated inclusions, 
and they send maximal 
(A2-maximal, qp-maximal) subsets into
 maximal (A2-maximal, qp-maximal) ones.
\end{proposition}

\begin{lemma}
\label{semisimpsat}
Let a semisimple group $G^s$ act on a 
factorial affine variety $Z$.
Then every good $G^s$-set $U \subseteq Z$ 
is $G^s$-saturated in $Z$.
\end{lemma}

\begin{proof}
Consider the quotient morphism 
$\pi \colon U \to U \quot G^s$,
and cover $U \quot G^s$ by affine
open subsets $V_i \subseteq U \quot G^s$.
Then each $U_i := \pi^{-1}(V_i)$ 
is affine, and hence, the complement 
$A_i := Z \setminus U_i$ is of pure 
codimension one in $Z$.
Since $Z$ is factorial affine,
$A_i$ is the set of zeroes of a function 
$f_i \in \Gamma(Z,\mathcal{O})$.

We claim that $f_i$ is $G^s$-invariant.
In fact, for any $z \in U_i$, the map 
$g \mapsto f_i(g \mal z)$ is an invertible
function on $G^s$, and thus, by semisimplicity
of $G^s$, is constant.
So, $U_i$ is the complement of the 
zero set of a $G^s$-invariant function,
and thus it is $G^s$-saturated in 
$Z$. Thus, $U$ as the union of the $U_i$,
is $G^s$-saturated as well.
\end{proof}

\begin{proof}[Proof of Proposition~\ref{thm:reduct}]
First, we check that the assignments are well 
defined.
Let $U \subseteq Z$ be a good $G$-set.
Then $U$ is as well a good $G^s$-set.
Lemma~\ref{semisimpsat} ensures that $U$ 
is $G^s$-saturated in $Z$.
Thus, $\pi(U)$ is open in $Y := Z \quot G^s$.
Moreover, the induced morphism 
$\pi(U) \to U \quot G$ is a good quotient
for the $T$-action, see 
Proposition~\ref{goodquotprops}.

If $V \subseteq Y$ is a good $T$-set, then 
$\pi^{-1}(V) \to V$ is a good quotient for 
the $G^s$-action, and thus
$\pi^{-1}(V) \to V \quot T$ is a good quotient
for the $G$-action,
see again Proposition~\ref{goodquotprops}.
Thus, the assignments are well defined.
Since every good $G^s$-set $U \subseteq Z$ is 
saturated with respect to $\pi \colon Z \to Y$,
they are moreover inverse to each other.

The fact that the assignments $U \mapsto \pi(U)$
and $V \mapsto \pi^{-1}(V)$ preserve saturated 
inclusion relies on the fact that we have the 
induced isomorphisms
$U \quot G \cong \pi(U) \quot T$ and 
$V \quot T \cong  \pi^{-1}(V) \quot G$. 
Moreover, this implies that maximality
(A2-maximality, qp-maximality) is preserved.
\end{proof}

\begin{proof}[Proof of Theorem~\ref{thm:maxcoll}]
In~\cite[Sec.~1]{ArHa}, the 
assertions were proven for torus actions 
on factorial affine varieties.
In particular, they hold 
for the action of $T := G/G^s$ 
on $Y := Z \quot G^s$.
Now, we have the canonical 
bijection between the respective sets
of orbit cones
$$ 
\Omega(Z) \ \to \ \Omega(Y), 
\qquad
\omega(z) \ \mapsto \ \omega(\pi(z)).
$$
In particular, this gives a one-to-one 
correspondence between the sets of
2-maximal collections $\Psi \subseteq \Omega(Z)$
and 
2-maximal collections $\Psi \subseteq \Omega(Y)$.
Thus, denoting by 
$V(\Psi) \subseteq Y$ the A2-maximal good $T$-set 
corresponding to a 2-maximal collection
 $\Psi \subseteq \Omega(Y)$,
Proposition~\ref{thm:reduct} reduces the 
problem to showing
$$ 
U(\Psi) \ = \ \pi^{-1}(V(\Psi)),
\qquad
\qquad
\Psi(U) \ = \ \Psi(\pi(U)).
$$
The first equality is obvious.
For ``$\subseteq$'' in the second one,
note that $U \mapsto \pi(U)$ is a good quotient
for the $G^s$-action.
Thus, if $G \mal z \subseteq U$ is closed, then 
$T \mal \pi(z) = \pi(G \mal z)$ is closed in 
$\pi(U)$.
For  ``$\supseteq$'', let $T \mal z \subseteq \pi(U)$
be closed. 
Then $\pi^{-1}(T \mal z)$
is a closed $G$-invariant subset in
$\pi^{-1}(\pi(U)) = U$ 
and thus contains a closed $G$-orbit, 
which is 
mapped onto $T \mal z$.
\end{proof}


\section{Lifting to the total coordinate space}
\label{sec:lift}

Here we reduce the problem of finding 
good $G$-sets for a given $G$-variety $X$ 
to the problem of finding good $(H \times G)$-sets,
$H$ a torus,
in a certain affine factorial variety $\b{X}$,
called the ``total coordinate space'' of $X$.
We begin with fixing the setup and recalling 
basic constructions from~\cite{BeHa1}.

Let $X$ be a normal algebraic variety 
with finitely generated free divisor 
class group $\Cl(X)$, and suppose that
$\Gamma(X,\mathcal{O}^*) = \KK^*$ holds.
To define the 
{\em Cox ring\/} (also {\em total coordinate ring\/})
$\mathcal{R}(X)$ of $X$,
choose a subgroup $K \subseteq \WDiv(X)$
of the group of Weil divisors mapping isomorphically 
onto $\Cl(X)$, and set
$$ 
\mathcal{R}(X)
\ := \
\Gamma(X,\mathcal{R}), 
\qquad
\text{where} 
\quad
\mathcal{R}
\ := \
\bigoplus_{D \in K} 
\mathcal{O}(D). 
$$
Then $\mathcal{R}(X)$ is a ring, where multiplication 
takes place in the field $\KK(X)$ of rational functions.
The definition of $\mathcal{R}(X)$ is 
(up to isomorphism) independent from 
the choice of  $K \subset \WDiv(X)$.
An important property of $\mathcal{R}(X)$ 
is that it admits unique factorization, 
compare~\cite{BeHa1}.

Throughout this section,
we assume that $\mathcal{R}(X)$ is 
finitely generated as a $\KK$-algebra;
this holds for spherical
varieties $X$, and, more generally 
for unirational
varieties $X$ with a complexity one group 
action, i.e., some Borel subgroup 
has an orbit of codimension one, 
see~\cite{Kn}.
We consider the following geometric
objects associated to the 
$K$-graded sheaf $\mathcal{R}$ of
$\mathcal{O}_X$-algebras:
$$
H \ := \ \Spec(\KK[K]),
\qquad
\b{X} \ := \ \Spec(\mathcal{R}(X)),
\qquad
\rq{X}
\ := \
\Spec_X(\mathcal{R}).
$$
The relative spectrum  $\rq{X}$ as well as 
$\b{X}$ come with actions of the Neron-Severi 
torus $H$,
both defined by the $K$-gradings of $\mathcal{R}$ 
and $\mathcal{R}(X)$ respectively.
The canonical morphism $q \colon \rq{X} \to X$
is a good quotient for the action of $H$,
and there is a canonical open $H$-equivariant 
embedding $\rq{X} \subseteq \b{X}$ 
with $\b{X} \setminus \rq{X}$ of codimension
at least two in $\b{X}$,
compare also~\cite{BeHa1}.
We will call $q \colon \rq{X} \to X$ a {\em Cox construction\/}
for $X$, as it naturally generalizes the 
often studied case of toric varieties~\cite{Co}.
Moreover, we refer to $\b{X}$ as the {\em total coordinate space}.

In this section, we consider $G$-equivariant
Cox constructions in the sense that a linear algebraic 
group $G$ acts on $\b{X}$ and $X$ such that the actions of 
$G$ and $H$ on $\b{X}$ commute,   
$\rq{X} \subseteq \b{X}$ is $G$-invariant and 
$q \colon \rq{X} \to X$ is $G$-equivariant.
The following two remarks can be helpful
for finding equivariant Cox constructions.

\begin{remark}
If a  connected 
linear algebraic group $G$ acts on $X$,
then the simply connected covering group
$\t{G}$ does as well.
After fixing a basis $E_1, \ldots, E_k$ 
of $K$ one may choose a $\t{G}$-linearization 
of each $E_i$. 
This induces a $\t{G}$-linearization of any 
$D \in K$ and thus defines a $\t{G}$-action on 
$\rq{X}$ making $\rq{X} \to X$ equivariant.
The lifted $\t{G}$-action extends to $\b{X}$.
Note that the actions of $G$ and $\t{G}$
on $X$ have the same quotients.
\end{remark}

\begin{remark}
Suppose that a linear algebraic group $G$ 
acts on a factorial affine variety $\b{X}$,
and that, moreover, there is an action of 
an algebraic torus $H$ on $\b{X}$ commuting
with the action of $G$.
Let $\rq{X} \subseteq \b{X}$ be invariant
under the actions of $H$ and $G$, 
and suppose that there is a good quotient 
$q \colon \rq{X} \to X$.
If there is an $H$-saturated subset 
$W \subseteq \rq{X} \subseteq \b{X}$
with $\b{X} \setminus W$ of codimension
at least two in $\b{X}$ such that
$H$ acts freely on $W$, 
then $q \colon \rq{X} \to X$ is an 
equivariant Cox construction for $X$.
\end{remark}

Given a $G$-equivariant Cox construction
$q \colon \rq{X} \to X$ 
with some reductive group $G$,
the (commuting) actions of $H$ and $G$
on $\b{X}$ define an action 
of the direct product
$H \times G$ on the factorial 
affine variety $\b{X}$.
Our aim is to relate the good 
$(H \times G)$-sets of $\b{X}$ 
to the good $G$-sets of $X$.
The key construction for this is 
the following.

\begin{definition}
\label{def:satint}
Let $G$ be a connected reductive group,
$X$ a $G$-variety with equivariant 
Cox construction 
$q \colon \rq{X} \to X = \rq{X} \quot H$
and total coordinate space $\b{X}$.
For every good $(H \times G)$-set
$W \subseteq \b{X}$, we set
\begin{eqnarray*}
W \sqcap_{G} \rq{X}
& := & 
\left\{
x \in W \cap \rq{X}; \;
\lim_{H \times G}(x,W) \subseteq \rq{X}, \
H \mal x_0 \text{ is closed} 
\right.
\\
& & 
\left.
\hphantom{\left\{ x \in W \cap \rq{X}; \; \right.} 
\text{in } \rq{X} \text{ for every }
x_0 \in \lim_{H \times G}(x,W) \right\}.
\end{eqnarray*}
\end{definition}

\begin{remark}
Let $G$ be a connected reductive group,
$X$ a $G$-variety with equivariant 
Cox construction 
$q \colon \rq{X} \to X = \rq{X} \quot H$
and total coordinate space $\b{X}$.
\begin{enumerate}
\item
If $X$ is $\QQ$-factorial, then $q \colon \rq{X} \to X$
is even a geometric quotient for the action of $H$ and,
for every good $(H \times G)$-set
$W \subseteq \b{X}$, one has
\begin{eqnarray*}
W \sqcap_{G} \rq{X}
& := & 
\left\{
x \in W \cap \rq{X}; \;
\lim_{H \times G}(x,W) \subseteq \rq{X} \right\}.
\end{eqnarray*}
\item
If $X$ is affine, then $\rq{X} = \b{X}$ holds and,
for every good $(H \times G)$-set
$W \subseteq \b{X}$, one has
\begin{eqnarray*}
W \sqcap_{G} \rq{X}
& := & 
\left\{
x \in W; \;
x_0 \in \lim_{H \times G}(x,W)
\Rightarrow 
H \mal x_0 \subseteq \b{X} 
\text{ closed}
 \right\}.
\end{eqnarray*}
\end{enumerate}
\end{remark}

The following result shows how to relate 
the good $(H \times G)$-sets of the total 
coordinate space $\b{X}$ 
to the good $G$-sets of $X$ using 
the assignment $W \mapsto W \sqcap_G \rq{X}$.

\begin{theorem}
\label{thm:sqcap}
Let $G$ be a connected reductive group,
$X$ a $G$-variety with equivariant 
Cox construction 
$q \colon \rq{X} \to X = \rq{X} \quot H$
and total coordinate space $\b{X}$.
Then, for every good $(H \times G)$-set 
$W \subseteq \b{X}$, the set 
$W \sqcap_G \rq{X}$ is $(H \times G)$-saturated 
in $W$ and $H$-saturated in $\rq{X}$.
This gives a surjection
\begin{eqnarray*}
\left\{
\text{good $(H \times G)$-sets of $\b{X}$}
\right\}
& \longrightarrow &
\left\{
\text{good $G$-sets of $X$}
\right\}
\\
W & \mapsto & q(W \sqcap_G \rq{X}).
\end{eqnarray*}
This map has $U \mapsto q^{-1}(U)$ as a right inverse.
Moreover, any maximal (A2-maximal, qp-maximal) 
good $G$-set $U \subseteq X$ 
is of the form $U = q(W \sqcap_G \rq{X})$
with a maximal (A2-maximal, qp-maximal) 
good $(H \times G)$-set $W \subseteq \b{X}$.
\end{theorem}

\begin{proof}
The first thing we have to show is 
that, for any good $(H \times G)$-set
$W \subseteq \b{X}$,
the set $W \sqcap_G \rq{X} \subseteq X$
is open and $(H \times G)$-saturated 
in $W$ and $H$-saturated in $\rq{X}$.
We do this by constructing 
$W \sqcap_G \rq{X} \subseteq W$
via stepwise removing suitable 
closed subsets from $W$. 

Let $p \colon W \to W \quot (H \times G)$ be the quotient,
and consider the closed $(H \times G)$-invariant subset
$A := W \setminus \rq{X}$ of $W$.
By the general properties of good quotients, 
we obtain an open, 
$(H \times G)$-saturated subset $V \subseteq W$ by
setting
$$ 
V 
\ := \
W \setminus p^{-1}(p(A))
\ = \  
\left\{
x \in W \cap \rq{X}; \; 
\lim_{H \times G}(x,W) \subseteq \rq{X}
\right\}
\ \subseteq  \ 
W \cap \rq{X}.
$$
Now, we consider the quotient 
$q \colon \rq{X} \to X$ and the 
$(H \times G)$-invariant, closed complement 
$B := \rq{X} \setminus V$.
Using $G$-equivariance and again the properties 
of good quotients, we obtain an $H$-saturated, 
$(H \times G)$-invariant open subset 
$V' \subseteq \rq{X}$ by setting 
$$ 
V' \ := \ 
\rq{X} \setminus q^{-1}(q(B))
\ = \ 
\left\{
x \in V; \; 
\lim_{H}(x,\rq{X}) \subseteq V
\right\}
\ \subseteq \
V.
$$
In the third shrinking step, we consider the quotient
$p \colon V \to V \quot (H \times G)$ and the 
closed $(H \times G)$-invariant subset $C := V \setminus V'$
of $V$.
Then we obtain an open $(H \times G)$-saturated subset
$V'' \subseteq V$ by setting 
$$ 
V''
\ := \ 
V \setminus p^{-1}(p(C))
\ = \ 
\left\{
x \in V; \; 
\lim_{H \times G}(x,V) \subseteq V'
\right\}
\ \subseteq \
V'.
$$
Thus, $V'' \subseteq W$ is $(H \times G)$-saturated
with good quotient 
$p \colon V'' \to V'' \quot (H \times G)$.
Moreover, $V'' \subseteq V'$ and hence 
$V'' \subseteq \rq{X}$
are $H$-saturated inclusions, because for 
$x \in V''$ the limit 
$\lim_H(x,V')$ is contained in the closure 
of $(H \times G) \mal x$ taken in $V$, which 
in turn is contained in $V''$.
In particular, we have a good quotient 
$q \colon V'' \to q(V'')$, where $q(V'') \subseteq X$ 
is open.
Hence, in order to finish the proof,
we have to verify
\begin{eqnarray*}
V'' & = & W \sqcap_G \rq{X}.
\end{eqnarray*}

Given $x \in V''$, we have 
$\lim_{H \times G}(x,W) = \lim_{H \times G}(x,V) \subseteq V'$.
In particular, $\lim_{H \times G}(x,W)$ is contained in $\rq{X}$.
Moreover, for $x_0 \in \lim_{H \times G}(x,W)$ one obtains
$\lim_H(x_0,\rq{X}) \subseteq V$, which gives 
$\lim_H(x_0,\rq{X}) \subseteq \lim_{H \times G}(x,V)$.
Since all $H$-orbits in $\lim_{H \times G}(x,V)$ 
have the same dimension in $\rq{X}$, we see that
$H \mal x_0$ is closed in $\rq{X}$.
Thus, $x \in W \sqcap_G \rq{X}$ holds.

Conversely, for any $x \in W \sqcap_G \rq{X}$, one 
obviously has $x \in V$.
Moreover, for $x_0 \in \lim_{H \times G}(x,W)$,
the orbit $H \mal x_0$ is closed in $\rq{X}$,
which gives $x_0 \in V'$.
This in turn means $x \in V''$.

Having seen that $W \sqcap_G \rq{X}$ is 
$(H \times G)$-saturated in $W$
and $H$-saturated in $\rq{X}$
for every good $(H \times G)$-set
$W \subseteq \b{X}$,
it is clear that we have the surjection
$W \mapsto q(W \sqcap_G \rq{X})$
as in the assertion.
Moreover, $U \mapsto q^{-1}(U)$ is obviously
a right inverse.

A few explaining words are needed concerning the claim 
that any maximal good $G$-set $U \subseteq X$ arises
as $U = q(W \sqcap_G \rq{X})$ with a maximal good
$(H \times G)$-set $W \subseteq \b{X}$.
In fact, $q^{-1}(U)$ is an $(H \times G)$-saturated 
subset of some maximal good
$(H \times G)$-set $W \subseteq \b{X}$.
Since $W \sqcap_G \rq{X} \subseteq W$ is 
$(H \times G)$-saturated  as well, we can
conclude that $q^{-1}(U) \subseteq W \sqcap_G \rq{X}$
is $(H \times G)$-saturated.
It follows that $U \subseteq q(W \sqcap_G \rq{X})$ 
is $G$-saturated, and thus, using maximality,
we obtain $U = q(W \sqcap_G \rq{X})$. 
\end{proof}

We now consider the 
setting of sets of semistable points.
This needs to recall a pullback 
construction for $G$-linearized divisors,
which was performed for the case of a torus $G$ 
in~\cite[Sec.~3]{BeHa2}
but generalizes without changes to 
any linear algebraic group $G$.

Let $D$ be any $G$-linearized Weil 
divisor on $X$.
Then the restriction $D_\reg$ to 
the set $X_\reg$ of regular points 
on $X$ is a linearized Cartier divisor,
and thus has a canonically 
$(H \times G)$-linearized 
pullback divisor $q^*D_\reg$,
where the $(H \times G)$-action on 
$$ 
q^{-1}(X_\reg)(q^*D_\reg)
\ \cong \ 
q^{-1}(X_\reg) \times_{X_\reg} X_\reg(D_\reg)
$$ 
is given by the diagonal $G$-action and 
the $H$-action on the first factor.
Since the complement $\b{X} \setminus q^{-1}(X_\reg)$
has codimension at least two, we may close the 
components of  $q^*D_\reg$, and obtain in this
way a $(H \times G)$-linearized Weil divisor $\b{D}$
on $\b{X}$. 
As shown in~\cite[Lemma~3.3]{BeHa2},
this construction sets up an isomorphism
\begin{eqnarray}
\label{eqn:pbiso} 
\Cl_G(X) \to \Cl_{H \times G}(\b{X}),
\qquad
[D] \ \mapsto \ [\b{D}].
\end{eqnarray}
The following statement shows that all 
sets of semistable points of $G$-linearized
divisors on $X$ arise from those of 
$(H \times G)$-linearized divisors on 
$\b{X}$; the proof is identical to that 
in the case of a torus $G$, 
see~\cite[Theorem~3.5]{BeHa2},
and therefore will be omitted.

\goodbreak

\begin{theorem}
\label{thm:sqcapsemistab} 
Let $G$ be a connected reductive group,
$X$ a $G$-variety with equivariant 
Cox construction
$q \colon \rq{X} \to X = \rq{X} \quot H$
and total coordinate space $\b{X}$.
Then, for any $G$-linearized Weil divisor
$D$ on $X$, we have a $(H \times G)$-saturated 
inclusion
$$
q^{-1}(X^{ss}(D,G))
\ = \
\b{X}^{ss}(\b{D}, H \times G) \sqcap_G \rq{X}
\ \subseteq \
\b{X}^{ss}(\b{D}, H \times G).
$$
\end{theorem}

An important finiteness result by Dolgachev and 
Hu~\cite{DoHu} and, independently, Thaddeus~\cite{Tha}
says that on any projective $G$-variety,
where $G$ is a reductive group, there are 
only finitely many GIT-equivalence classes arising from 
{\em ample\/} bundles.
In our setting, Theorem~\ref{thm:sqcapsemistab} 
gives more:

\begin{corollary}
Let $G$ be a connected reductive group, and 
$X$ a $G$-variety with finitely generated total
coordinate ring.
Then the $G$-action on $X$ has only finitely many 
GIT-equivalence classes.
\end{corollary}

\begin{proof}
According to Theorem~\ref{thm:sqcapsemistab},
the number of GIT-equivalence classes of the $G$-action on $X$ 
is boundend by the number of GIT-classes of 
the $(H \times G)$-action on $\b{X}$.
But the latter number is finite by 
Theorem~\ref{thm:GITfan}.
\end{proof}

\section{Computing a first example}
\label{sec:firstex}

The previous sections suggest the following
strategy for constructing good $G$-sets
of a given $G$-variety $X$.
First, take an equivariant Cox 
construction 
$q \colon \rq{X} \to X = \rq{X} \quot H$ 
and consider the associated 
total coordinate space $\b{X}$.
Then Theorem~\ref{thm:sqcap} reduces the 
problem of finding good $G$-sets 
$U \subseteq X$ to finding good 
$(H \times G)$-sets of $\b{X}$.
By Proposition~\ref{thm:reduct}, the latter 
problem is equivalent to finding the 
quotients of a torus action on 
$\b{Y} = \b{X} \quot G^s$.

In fact, in many concrete cases, 
the equivariant Cox construction is
given from the beginning, and 
Classical Invariant Theory often
provides enough information
on the quotient $\b{X} \quot G$,
see the examples treated later.
So the general difficulties remain in 
understanding the step 
$W \mapsto W \sqcap_G \rq{X}$ 
and the computation of 
GIT-fan and A2-maximal collections
for torus actions.

The first problem dissappears,
for example,
when we restrict to GIT-quotients
arising from ample bundles,
see Section~\ref{sec:amplegit}.
For the second one, we begin
with a general observation showing that 
one may work in terms 
of {\em walls\/}, i.e., orbit cones 
of codimension one.
Let us first spend a few words on 
the combinatorial framework. 

\begin{remark}
\label{rem:wallgen}
Let $\omega_1, \ldots, \omega_r$ be
polyhedral cones 
in a rational vector space $K_\QQ$
such that their union is
a convex cone $\omega \subseteq K_\QQ$.
Suppose that
$$ 
\Sigma 
\ := \ 
\Sigma(\omega_1, \ldots, \omega_r)
\ := \
\{\lambda(u); \; \ u \in \omega\},
\qquad
\text{where }
\lambda(u) \ := \ \bigcap_{u \in \omega_i} \omega_i,
$$
is a fan, any $\omega_i$ is a face of some full 
dimensional $\omega_j \subseteq K_\QQ$, and 
the facets $\eta_1, \ldots, \eta_s$ 
of the full dimensional $\omega_j$ occur among the 
$\omega_i$. Then
\begin{itemize}
\item
the maximal cones $\lambda(u) \in \Sigma$ 
are precisely the closures of the connected components of 
$\omega \setminus (\eta_1 \cup \ldots \cup \eta_s)$,
\item
every nonmaximal cone $\lambda(u) \in \Sigma$ is the intersection 
over the facets $\eta_i$ with $\lambda(u) \subseteq \eta_i$. 
\end{itemize}
\end{remark}

If we are in the setting~\ref{rem:wallgen}, 
then we call the facets $\eta_1, \ldots, \eta_s$
of the full-dimensional $\omega_j$ the {\em walls\/} 
and we say that $\Sigma$ is {\em determined by the walls}.

\begin{proposition} 
\label{lem:walls}
Let a reductive group $G$ act on a factorial affine
variety $Z$.
Then the associated GIT-fan is determined by 
its walls.
\end{proposition}

\begin{proof}
According to Lemma~\ref{lem:reduct}, 
we may assume that
$G$ is a torus, acting effectively.
To obtain the setting~\ref{rem:wallgen},
two things have to be verified.
Firstly, given an orbit cone of full dimension,
then also its facets are orbit cones;
this is obvious.
Secondly, every orbit cone is a face of some 
orbit cone of full dimension;
this will be done below.

We have to show that any $G$-orbit is contained 
in the closure of a $G$-orbit of maximal dimension. 
Otherwise, we find some orbit $G \mal z$
such that $\dim (G \mal z)$ is not maximal 
and $G \mal z$ is not contained in the closure 
of any other $G$-orbit. 
Then $G' := (G_z)^0$ is a proper subtorus of 
$G$, and $G'$ acts nontrivially on $Z$.
Semicontinuity of fiber dimension tells us
that the fiber $\pi^{-1}(\pi(z))$ of the quotient map
$\pi \colon Z \to Z \quot G'$ must contain a 
$G'$-orbit $G' \mal z'$ of positive dimension.
As a $G'$-fixed point, $z$ lies in the closure
of $G' \mal z'$.
It follows that $G \mal z$ is contained in the 
closure of the orbit $G \mal z'$, which is 
different from $G \mal z$; a contradiction.
\end{proof}

\begin{example}
Consider the homogeneous space $X := \SL(3) / H$,
where $H \subset \SL(3)$ is a maximal torus.
Then $X$ is a smooth affine variety of dimension~6,
and the special orthogonal group $G := \SO(3) \subset \SL(3)$
acts on $X$ from the left. 
The generic $G$-orbit on $X$ is of dimension 3 and it 
is closed in $X$, see~\cite{Lu}.
  
Consider $\rq{X} := \SL(3)$ with the 
left $G$-action.
Then the projection 
$q \colon \rq{X} \to X$ is a $G$-equivariant Cox 
construction, and the total coordinate space is
given as $\b{X} = \rq{X}$.
Moreover, $\b{Y} := \b{X} \quot G$ is the
homogeneous space $\SO(3)\backslash \SL(3)$ 
with respect to the left $\SO(3)$-action.
The situation is summarized in the following 
commutative diagram
$$
\xymatrix{
& 
\b{X}=\hat X =\SL(3) \ar[dr]^{\pi} \ar[dl]_{q} & 
\\
X=\SL(3)/H \ar[dr]_{\varphi}
& & 
\b{Y} = \SO(3)\backslash \SL(3) \ar[dl]^{\psi} 
\\
& Y := \SO(3)\backslash \SL(3)/H.}
$$

Combining Proposition~\ref{thm:reduct} and Theorem~\ref{thm:sqcap}, 
we see that, in the present setting,
the good $G$-sets of $X$ are in bijection with 
the good $H$-sets of $\b{Y}$ via 
$U \mapsto \pi(q^{-1}(U))$. 
Moreover, the quotient of $U$ by $G$ is geometric 
if and only if the quotient of 
$\pi(q^{-1}(U))$ by $H$ is so.
So, our task is reduced to describing the 
$H$-quotients of $\b{Y}$.

First recall that $\b{Y}$ can be identified 
as the variety of symmetric 
$(3\times 3)$-matrices with determinant one
via 
$$G \mal A \ \mapsto \ A^t \cdot A.$$
The $H$-action is given as 
$(t_1,t_2,t_3)(a_{ij})=(t_it_ja_{ij})$,
where $t_3 = t_1^{-1}t_2^{-1}$.
The following matrices 
have one-dimensional $H$-orbits:
$$
A_1 \ := \ 
\left(
\begin{array}{rrr}
  -1 & 0 & 0 \\
  0 & 0 & 1 \\
  0 & 1 & 0
\end{array}
\right), 
\qquad
A_2 \ := \ 
\left(
\begin{array}{rrr}
  0 & 0 & 1 \\
  0 & -1 & 0 \\
  1 & 0 & 0
\end{array}
\right).
\qquad
A_3 \ := \ 
\left(
\begin{array}{rrr}
  0 & 1 & 0 \\
  1 & 0 & 0 \\
  0 & 0 & -1
\end{array}
\right).
$$
The associated orbit cones in $\XX(H) = \QQ^2$ 
are the lines $\omega(A_1) = \QQ \mal e_1$,
$\omega(A_2) = \QQ \mal e_2$ and
$\omega(A_3) = \QQ \mal (e_1+e_2)$.
According to Proposition~\ref{lem:walls} 
the GIT-fan looks as follows.
\begin{center}
\input{butterfly1.pstex_t}
\end{center}
The six full dimensional cones of this fan define
geometric GIT-quotients 
$U_i \to U_i / G$,
where $1 \le i \le 6$,
of the $G$-variety
$X$. 
Moreover, there are two 2-maximal collections
defining A2-maximal quotients 
$U_i \to U_i / G$, where $7 \le i \le 8$,
with a non-quasiprojective 
quotient space, namely
$$ 
(\KK_\QQ,\omega_1^+,\omega_2^+,\omega_3^+),
\qquad \qquad
(\KK_\QQ,\omega_1^-,\omega_2^-,\omega_3^-),
$$
where $\omega_i^+$ are half spaces
bounded by $\omega(A_i)$, and $\omega_i^- = -\omega_i^+$
such that the intersection over the interiors of 
the $\omega_i^+$ as well as that over the interiors
of the $\omega_i^-$ are empty.
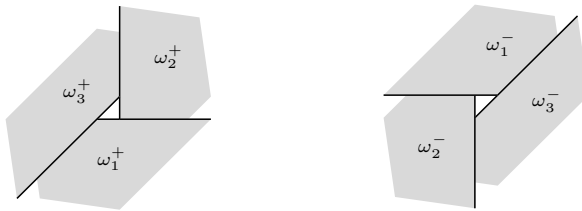
\begin{center}
\input{butterfly3.pstex_t}
\qquad
\qquad
\qquad
\input{butterfly2.pstex_t}
\end{center}

It is not hard to check that the quotient
$X \to X \quot G$ has precisely three exceptional
fibers, each of which consists of one closed orbit
and two one-parameter families of three-dimensional
orbits.
From this one may guess that there 
are eight maximal open subsets with geometric 
quotient. 
In any case, the quotient space $U_i/G$ is a small modification
of the affine threefold $X \quot G$, having three
exceptional fibers, each isomorphic to a projective line.
A priori it is not clear why six of the $U_i/G$ should be 
quasiprojective and two not.
\end{example}

\section{Chambers of the linearized ample cone}
\label{sec:amplegit}

In this section, we consider projective 
$G$-varieties.
In~\cite{DoHu} and~\cite{Tha} it was 
first proved that the cone of linearized
ample divisors has a ``chamber structure'' 
describing the GIT-equivalence, see also
\cite{Re}.
In this section, we describe this chamber 
structure in our setting and then treat
a concrete example.

Let $G$ be a connected reductive group,
and let $X$ be a projective $G$-variety with equivariant
Cox construction $q \colon \rq{X} \to X = \rq{X} \quot H$ 
and total coordinate space $\b{X}$.
Let $K$ denote the character group of the torus $H$
and $M$ that of $G$. 
Then we have canonical isomorphisms
$$ 
\Cl(X) \ \cong \ K,
\qquad \qquad
\Cl_G(X) 
\ \cong \ 
\Cl_{H \times G}(\b{X})
\ \cong \
K \times M,
$$ 
see formulae~(\ref{eqn:lin2char}) and~(\ref{eqn:pbiso}) 
for the latter two.
This allows us to denote the $G$-linearized divisors 
on $X$ by pairs $(D,\chi) \in K \times M$.
Moreover, we denote by $\kappa_X \subseteq K_\QQ$ 
the cone of semiample divisor classes on $X$; recall
from~\cite{BeHa2} 
that $\kappa_X$ is a GIT-cone for the $H$-action on 
$\b{X}$. Finally, we denote by 
$\omega(\b{X}) \subseteq K_\QQ \times M_\QQ$ 
the weight cone and by $\Sigma(\b{X})$ the GIT-fan 
of the $(H \times G)$-action 
on $\b{X}$.

\begin{proposition}
\label{amplegit}
Let $G$ be a connected reductive group,
and $X$ a projective
$G$-variety with equivariant
Cox construction $q \colon \rq{X} \to X = \rq{X} \quot H$ 
and total coordinate space $\b{X}$.
Then, in the above notation, the following holds.
\begin{enumerate}
\item
The cone $\alpha(X) \subseteq K \times M$ 
of ample $G$-linearized 
divisor classes with nonempty set of semistable 
points is given by
\begin{eqnarray*}
\alpha(X) 
& = & 
(\kappa_X^{\circ} \times M_\QQ) \ \cap \ \omega(\b{X}).
\end{eqnarray*}
\item
The partial fan
$\Sigma(X) := 
\{\lambda \cap \alpha(X); \lambda \in \Sigma(\b{X}) \}$
describes the GIT-equivalence on $X$
in the sense that for any two
$(D,\chi)$ and $(D',\chi')$ in $\alpha(X)$,
one has
\begin{eqnarray*}
X^{ss}(D,\chi) \subseteq X^{ss}(D',\chi')
& \iff & 
\lambda(D,\chi) \succeq \lambda(D',\chi').
\end{eqnarray*}
\end{enumerate}
\end{proposition}

\begin{proof}
We prove~(i). 
Given any $G$-linearized divisor 
on $X$, represented by some 
$(D,\chi) \in K \times M$, 
its invariant sections are exactly the 
semi-invariants with respect to the
weight $(D,\chi)$ in $\Gamma(\b{X},\mathcal{O})$.
So, the weight cone $\omega(\b{X})$ contains 
precisely the $G$-linearized divisors of $X$
admitting invariant sections.

In order to verify the description of 
$\alpha(X)$, it suffices to show that 
for any $(D,\chi) \in \omega(\b{X})$ 
with $D$ ample, some positive multiple 
$nD$ admits a section $f$ with 
$X \setminus Z(f)$ affine.
But this follows from the general 
observation that for any section 
of an ample divisor the complement 
of its zero set is affine.

To see assertion~(ii), note that
for any ample $G$-linearized divisor $D$
on the projective variety $X$, 
the quotient space $X^{ss}(D) \quot G$ 
is again projective.
Thus Theorem~\ref{thm:sqcapsemistab}
gives
$q^{-1}(X^{ss}(D)) = \b{X}^{ss}(\b{D}, H \times G)$. 
Consequently, the assertion follows from 
Theorem~\ref{thm:GITfan}.
\end{proof}

We now treat a concrete example.
For $n \ge 2$, let $G=\Sp(2n)$ be the symplectic group, 
i.e., the group of invertible matrices preserving
a non-degenerate skew-symmetric bilinear form $\bangle{\ , \ }$ 
on $\KK^{2n}$. 
Then $G$ acts diagonally on the $m$-fold product
$$
X \ := \ (\PP^{2n-1})^m.
$$ 
Fix a hyperplane $E_i$ on each factor $\PP^{2n-1}$ and 
consider its pullback $D_i$ on $X$.
Then the lattice $K \subseteq \WDiv(X)$ 
generated by $D_1, \ldots, D_m$ maps isomorphically 
to the divisor class group $\Cl(X)$.
We have the identification  
$K \cong \ZZ^m$ via $D_i \mapsto e_i$ and,
moreover,
$$
\b{X} 
\ = \ 
\Spec\left( \bigoplus_{D \in K} \Gamma(X,\mathcal{O}(D))\right)
\ \cong \
(\KK^{2n})^m
\ =: \ 
V.
$$
The torus acting on $\b{X}$ is $H \cong (\KK^*)^m$,
and its action is componentwise scalar multiplication. 
As $G$ has trivial character group, the GIT-fan
of the $(H \times G)$-action lives in 
$K_\QQ = \XX_\QQ(H) \cong \QQ^m$.

For our description of the GIT-fan,
we need one more notation.
Given a set of vectors 
$w_1, \ldots, w_k$,
consider the collection 
$\omega_1, \ldots, \omega_r$ of all
convex cones generated by some of the $w_i$,
and set 
\begin{eqnarray*}
\Sigma(w_1, \ldots, w_k) 
& := & 
\Sigma(\omega_1, \ldots, \omega_r)
\end{eqnarray*}
Note that $\Sigma(w_1, \ldots, w_k)$
is the coarsest common refinement of all fans 
having precisely 
$\QQ_{\ge 0} \mal w_1, \ldots, \QQ_{\ge 0} \mal w_k$
as their set of rays.

\begin{theorem}
\label{tsp}
For $1 \le i < j \le m$, let 
$u_{ij} = (u_{ij}^1, \ldots, u_{ij}^m) \in \ZZ^m$ 
be the vector with entries 
$1$ at the $i$-th and $j$-th
place and $0$ elsewhere.
Then the weight cone $\omega(\b{X})$ 
of the $(H \times G)$-action  
on $\b{X}$ and the 
$G$-ample cone $\alpha(X)$ of the 
$G$-action on $X$ are 
\begin{eqnarray*}
\omega(\b{X}) 
& = &
\cone(u_{ij}; \; 1 \le i < j \le m)
\\
& = & 
\left\{
(s_1, \ldots, s_m) \in \QQ_{\ge 0}^m; \; 
2 s_i \le s_1 + \ldots + s_m, \ 1 \le i \le m
\right\}
\\[2ex]
\alpha(X)
& = & 
\QQ_{>0}^m \ \cap \ \omega(\b{X}) .
\end{eqnarray*}
The GIT-fan of the action of $H \times G$ on $\b{X}$  
is $\Sigma(u_{ij}; \; 1 \le i < j \le m)$;
it is determined by the walls, and these are 
precisely the orbit cones
\begin{eqnarray*}
 \omega(\b{X}) 
\cap \left\{(s_1,\dots,s_m); \; \sum_{j\in J}s_j=\sum_{l\notin J}s_l \right\},
& & 
J\subset\{1,\dots,m\},
\
1<|J|\le\frac{m}{2},
\\
\cone\left( u_{ij}; \; 
\sum_{k\in J_1}u^k_{ij} 
=
\sum_{l\in J_2}u^l_{ij} \right),
& & 
J_1,J_2\subset\{1,\dots,m\}, \
J_1\cap J_2=\emptyset, 
\\ 
& & 
|J_1|+|J_2|\le m-3.
\end{eqnarray*} 
\end{theorem}


\begin{proof}
Lemma~\ref{lem:reduct} tells us that the 
GIT-fan of the action of $H \times G$ on $\b{X}$
is the same as that of the $H$-action on 
$\b{Y} = \Spec(\KK[V]^G)$.
The algebra of invariants 
$\KK[V]^G$ is generated by the functions 
$f_{ij} \in \KK[V]$ given by 
\begin{eqnarray*}
f_{ij}(v_1,\ldots,v_m)
&  := & 
\bangle{v_i,v_j}.
\end{eqnarray*}
Each $f_{ij}$ is $H$-semiinvariant with 
weight $u_{ij}$.
Moreover, $\b{Y} = \Spec (\KK[V]^G)$ is the 
variety of skew-symmetric 
$m\times m$-matrices of rank $\le 2n$, 
see~\cite[Sec.~9]{PV2}. 
The quotient morphism
$\pi \colon V \to \b{Y}$ sends $(v_1,\dots,v_m)$ to 
the matrix $(\bangle{v_i,v_j})$, 
and an element 
$(t_1,\dots,t_m)\in H$
moves $(c_{ij})$ 
to $(t_it_jc_{ij})$.

According to Proposition~\ref{lem:walls}, the problem
of describing the GIT-fan is reduced to 
computing the walls of the $H$-action $\b{Y}$. 
To any subset $J \subset\{1,\dots,m\}$ 
with $1<|J| \le m/2$, we associate
an hyperplane
\begin{eqnarray*}
\mathcal{H}_J
& := &
\left\{(s_1,\dots,s_m); \; \sum_{j\in J}s_j=\sum_{l\notin J}s_l \right\}.
\end{eqnarray*}
Moreover, for any pair 
$J_1,J_2 \subset \{1,\dots,m\}$ 
of disjoint subsets satisfying
$|J_1|+|J_2|\le m-3$
and if $J_1$ is empty then $J_2 = \{i\}$, 
we set 
\begin{eqnarray*}
\mathcal{H}_{J_1,J_2}
&  := &
\left\{(s_1,\dots,s_m); 
\sum_{j_1 \in J_1}s_{j_1} 
=
\sum_{j_2 \in J_2} s_{j_2} 
\right\}.
\end{eqnarray*}

With these definitions, the description 
of the walls, i.e., the orbit cones of codimension
one, is an immediate consequence of the following 
two claims.

\medskip

\noindent {\em Claim 1.~}
If $C$ is a skew-symmetric matrix with a one-dimensional stabilizer $H_C$, 
then $\omega(C)$ lies in either some 
$\mathcal{H}_J$ or some $\mathcal{H}_{J_1,J_2}$.

\medskip

Let $C=(c_{ij})$.
First observe, that if $(t_1,\dots,t_m)$ stabilizes $C$, 
then for any two $i,j$ with $c_{ij}\ne 0$
we have $t_i=t_j^{-1}$. 
Next we associate a graph $\Gamma_C$ to $C$:
the set of vertices is $\{1,\dots,m\}$,
and the edges are the $(ij)$ 
with $c_{ij}\ne 0$.
Let $\Gamma_C^1, \ldots, \Gamma_C^k$ be the 
connected components of $\Gamma_C$. 

If $\Gamma_C^s$ contains a cycle of odd
length (type I), 
then $t^2_i=1$ holds for all vertices 
$i$ in $\Gamma_C^s$. 
If $\Gamma_C^s$ contains no cycle of odd
length (type II), then one may divide the 
set of vertices of $\Gamma_C^s$ into 
subsets $J_1,J_2$ with $J_1\cap J_2=\emptyset$ 
such that for any edge $(ij)$  
of $\Gamma_C^s$ we have $i\in J_1$,
$j\in J_2$.

Any connected component of type II gives a free
parameter in $H_C$. Thus, if the stabilizer 
$H_C$ is one-dimensional, there is exactly one 
connected component, say $\Gamma_C^1$, of type II 
and all others are of  type I. 
If $\Gamma_C^1=\Gamma_C$, then we have  
$\omega(C)\subset\mathcal{H}_{J_1}=\mathcal{H}_{J_2}$, 
and otherwise we have 
$\omega(C)\subset\mathcal{H}_{J_1,J_2}$
(any component of type I contains $\ge 3$ vertices).
This proves Claim 1.

\medskip

\noindent{\em Claim 2.~}
For any $\mathcal{H}_J$ 
(resp. $\mathcal{H}_{J_1,J_2}$),
there exists a skew-symmetric 
matrix $C$ of rank $\le 4$ such that
$\omega(C)$ is generated by all 
$u_{ij}\in\mathcal{H}_J$ 
(resp. $u_{ij}\in\mathcal{H}_{J_1,J_2}$), 
and $\omega(C)$ generates
$\mathcal{H}_J$ (resp. $\mathcal{H}_{J_1,J_2}$).

\medskip 

First consider  $\mathcal{H}_J$.
By renumbering, we may assume that 
$J=\{1,2,\ldots,k\}$. 
Then the hyperplane $\mathcal{H}_J$ 
is generated by $\omega(C(k,m-k))$,
where we set 
\begin{eqnarray*}
C(k,l)
& := & 
\left(
\begin{array}{cc}
0 & \mathds{1}_{k\times l}
\\
-\mathds{1}_{l \times k} & 0
\end{array}
\right)
\end{eqnarray*}
and $\mathds{1}_{k\times l}$ denotes 
the $(k\times l)$-matrix with all 
entries equal one.
One easily sees that all weights
$u_{ij}$ lying in the hyperplane $\mathcal{H}_J$ 
already belong to  $\omega(C(k,m-k))$.

Now we turn to $\mathcal{H}_{J_1,J_2}$.
Again, by renumbering, we may assume that
$J_1=\{1,\dots,k_1\}$ and 
$J_2=\{k_1+1, \ldots, k_1 + k_2\}$ holds.
Set $s := m-k_1-k_2$, and take 
pairwise non-proportional 
vectors $w_1,\dots,w_s$ 
in some twodimensional symplectic vector space 
$W$; then these vectors define a skew symmetric
matrix  $C(s)  = (\bangle{w_i,w_j})$ 
of rank two having only non-zero 
non-diagonal elements. 
The hyperplane $\mathcal{H}_{J_1,J_2}$ is generated by 
$\omega(C(J_1,J_2))$, where
\begin{eqnarray*}
C(J_1,J_2)
& = & 
\left(
\begin{array}{cc}
C(k_1,k_2) &  0 
\\
0   & C(s)
\end{array}
\right).
\end{eqnarray*}
Again one directly checks 
that all weights $u_{ij}$ in $\mathcal{H}_{J_1,J_2}$,
are already contained in the orbit cone 
$\omega(C(J_1,J_2))$.
This proves Claim 2.

Finally, observe that for $n\ge 2$ the GIT-fan of the action 
$G = \Sp(2n)$ on $(\PP^{2n-1})^m$ does not 
depend on $n$. 
But for $2n\ge m$, the variety $\b{Y}$ is a vector space, 
and hence, the corresponding GIT-fan coincides with
$\Sigma(u_{ij}; \; 1 \le i,j \le m)$. 
\end{proof}


\begin{remark}
It is proved in~\cite[Prop.~17]{Re} 
that for an $\SL(2)$-action on a projective
variety any wall is the intersection of 
the weight cone with a hyperplane.
In the setting of Theorem~\ref{tsp} 
shows that for $G=\Sp(2n)$, $n\ge 2$,
this is not the case. 
Indeed, the intersection 
$\mathcal{H}_{J_1,J_2}\cap\omega(\b{X})$
has extremal rays different from any 
of the $u_{ij}$; e.g. the ones generated by 
$$
 (\underbrace{0,\dots,0,1,0,\dots,0}_{|J_1|},
\underbrace{0,\dots,0,1,0,\dots,0}_{|J_2|},0,\dots,0,2,0,\dots,0).
$$
\end{remark}

\begin{remark}
In the setting of Theorem~\ref{tsp}, none of the quotients
$X^{ss}(D)\to X^{ss}(D)\quot G$ is geometric.
\end{remark}


A further class of examples arises from 
(reducible) representations of simple 
groups having a free algebra of invariants.
They are all known and can be found
in~\cite{AG} and \cite{Sch};
the multidegrees of basic invariants 
are also indicated in tables there.

\begin{example}
Consider the action of the special 
linear group $G := \SL(6)$ on
the product
$X = \PP(\KK^6) \times \PP(\Lambda^2\KK^6) \times \PP(\Lambda^3\KK^6)$.
Then the total coordinate space is 
$\b{X} = \KK^6 \times \Lambda^2\KK^6 \times \Lambda^3\KK^6$,
and $H = (\KK^*)^3$ acts by scalar multiplication 
on the factors.
The weights of the canonical generators of the algebra $\KK[\b{X}]^G$ 
are listed as well in \cite{Sch}; they are 
$$ 
w_1=(0,0,4),
\quad
w_2=(0,3,0),
\quad
w_3=(0,3,4),
$$
$$
w_4=(1,1,1),
\quad
w_5=(2,2,2),
\quad
w_6=(1,1,3).
$$
As the algebra  $\KK[\b{X}]^G$ of invariants 
is a polynomial ring,
the GIT-fan of the $H$-action on $\b{Y} = \b{X} \quot G$ 
is $\Sigma(w_1, \ldots, w_6)$. Here comes a figure,
showing the intersection of the weight cone $\omega(\b{Y})$ 
with a transversal hyperplane.
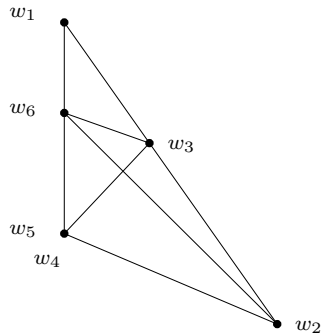
\begin{center}
\input{schwarzex.pstex_t}
\end{center}
All resulting quotients are toric varieties; 
the computation of the respective fans is 
a standard calculation.
\end{example}

\section{Gelfand-MacPherson type correspondences}
\label{sec:gmcorr}

The classical Gelfand-MacPherson correspondence~\cite{GM}
relates generic orbits of the diagonal
action of $\SL(n)$ on $(\PP^{n-1})^m$
to generic orbits of a torus action on the 
Grassmannian $G(n,m)$.
This may even be extended to isomorphisms
between certain quotient spaces on both
sides, see~\cite[Theorem~2.4.7]{Ka}.
Combining Proposition~\ref{thm:reduct} 
and Theorem~\ref{thm:sqcap}, 
we obtain the following general way to relate 
quotients for a reductive group action
to quotients of a torus action.

\begin{construction}
\label{rem:gelmac}
Let $G_X$ be a connected reductive group,
$X$ a $G_X$-variety with equivariant 
Cox construction
$q_X \colon \rq{X} \to X = \rq{X} \quot H_X$
and total coordinate space $\b{X}$.
Consider the induced $T$-action on $\b{Y}$,
where
$$
T 
\ := \ 
(H_X \times G_X)/(H_X \times G_X)^s 
\ = \ 
H_X \times (G_X/G_X^s),
$$
$$
\b{Y} 
\ := \ 
\b{X} \quot (H_X \times G_X)^s 
\ = \
\b{X} \quot G_X^s.
$$
Suppose that for some $T$-invariant open 
set $\rq{Y} \subseteq \b{Y}$ and some subtorus
$H_Y \subseteq T$ we obtain a Cox construction 
$q_Y \colon \rq{Y} \to Y = \rq{Y} \quot H_Y$,
and consider the induced action of 
$T_Y := T/H_Y$ on $Y$.
Then the good $G_X$-sets $U \subseteq X$ 
and the good $T_Y$-sets $V \subseteq Y$ 
fit into the diagram
$$ 
\xymatrix{
{\rq{U}}
\ar@{}[r]|{\subseteq}
\ar[d]_{\quot H_X}
&
{\rq{X}}
\ar@{}[r]|{\subseteq}
\ar[d]^{\quot H_X}
&
{\b{X}}
\ar[rr]^{\pi}_{\quot G_X^s}
& &
{\b{Y}}
&
{\rq{Y}}
\ar@{}[l]|{\supseteq}
\ar[d]_{\quot H_Y}
&
{\rq{V}}
\ar@{}[l]|{\supseteq}
\ar[d]^{\quot H_Y}
\\
{{U}}
\ar@{}[r]|{\subseteq}
\ar[d]
&
{{X}}
&
& &
&
{{Y}}
&
{{V}}
\ar@{}[l]|{\supseteq}
\ar[d]
\\
{{U \quot G_X}}
&
&
& &
&
&
{{V \quot T_Y}}
}
$$
where we set $\rq{U} := q_X^{-1}(U)$ and 
$\rq{V} := q_Y^{-1}(V)$.
Moreover, combining Proposition~\ref{thm:reduct} 
and Theorem~\ref{thm:sqcap} 
gives canonical assignments from good $G_X$-sets
$U \subseteq X$
to good $T_Y$-sets 
$V \subseteq Y$
and vice versa:

\begin{enumerate}
\item
If $U \subseteq X$ is a good $G_X$-set,
then 
$V := q_Y(\pi(\rq{U}) \sqcap_{T_Y} \rq{Y})$
is a good $T_Y$-set in $Y$, 
and there is a canonical open embedding 
$V \quot T_Y \to U \quot G_X$.
This embedding is an isomorphism if and only if 
one has an $H_Y$-saturated inclusion
$\pi(\rq{U}) \subseteq \rq{Y}$.
\item
If $V \subseteq Y$ is a good $T_Y$-set,
then 
$U := q_X(\pi^{-1}(\rq{V}) \sqcap_{G_X} \rq{X})$ 
is a good $G_X$-set in $X$, and 
there is a canonical open embedding
$U \quot G_X  \to V \quot T_Y$.
This embedding is an isomorphism if and only if 
one has an $H_X$-saturated inclusion
$\pi^{-1}(\rq{V}) \subseteq \rq{X}$.
\end{enumerate}
To consider these assignments for sets of
semistable points recall first
that~(\ref{eqn:lin2char}) and~(\ref{eqn:pbiso}) 
provide canonical isomorphisms
relating the respective groups of linearized 
divisor classes
$$
\Cl_{G_X}(X) 
\ \cong \ 
\Cl_{H \times G_X}(\b{X})
\ \cong \
\Cl_{T}(\b{Y}) 
\ \cong \
\Cl_{T_Y}(Y).
$$
If $U \subseteq X$ is a set of semistable points
of a $G_X$-linearized divisor, then Lemma~\ref{lem:reduct}
and Theorem~\ref{thm:sqcapsemistab} ensure that
the associated set $V \subseteq Y$ is a saturated 
(possibly proper) subset of the set of semistable 
points of the corresponding $T_Y$-linearized divisor 
and vice versa.
\end{construction}

In certain cases, the above Gelfand-MacPherson 
type correspondence can even be extended to 
the respective inverse limits of the GIT-quotients,
which in turn gives interesting descriptions of 
moduli spaces, see~\cite{Tha2}.
For giving a general statement in this context,
we now recall the basic facts on inverse limits 
of GIT-quotients.

Consider a projective variety $X$ with an action of 
a connected reductive group $G$. 
If, for two ample $G$-linearized divisors $D,D'$, 
we have 
an inclusion $X^{ss}(D) \subseteq X^{ss}(D')$, then 
there is an induced morphism of the 
associated quotient spaces.
These induced morphisms of GIT-quotients 
form an inverse system, 
the {\em ample GIT-system\/} of the $G$-variety $X$.
The {\em GIT-limit\/} of $X$ is the inverse limit of 
this system. 
As in the case of fiber products,
the GIT-limit can be realized as a subvariety 
of the product over all GIT-quotients arising from ample 
bundles.

In order to compare GIT-limits in our setting,
recall that for a projective $G$-variety $X$ 
with equivariant Cox construction 
$\rq{X} \to X = \rq{X} \quot H_X$
its ample cone is 
the relative interior of a GIT-cone $\kappa_X$ of the 
$H_X$-action on the total coordinate space
$\b{X}$. 
By the {\em open $G$-ample cone\/}, we mean 
the relative interior of the cone
$(\kappa_X \times M_\QQ) \cap \omega(\b{X})$ in 
$(K \times M)_{\QQ}$, where
$K$ and $M$ stand for the character lattices 
of $H_X$ and $G$ respectively.

\begin{theorem}
Consider a $G_X$-variety $X$ and a 
$T_Y$-variety $Y$ as in~\ref{rem:gelmac},
and suppose that both are projective.
If the canonical isomorphism
$\Cl_{G_X}(X) \to \Cl_{T_Y}(Y)$
sends the open $G_X$-ample cone 
onto the open $T_Y$-ample cone,
then the GIT-limits of $X$ and $Y$ 
are isomorphic.
\end{theorem}

\begin{proof}
First note that for determining the GIT-limit,
it suffices to consider GIT-quotients given by 
the classes inside the open linearized ample cone.
Any class inside the open $G_X$-ample cone defines 
an ample class on $Y$ and vice versa.
Moreover, since $X$ and $Y$ are projective, all 
quotients arising from ample classes are projective
again. This implies 
$$ 
q_X^{-1}(X^{ss}(D,G_X)) = \b{X}^{ss}(\b{D}, H_X \times G_X),
\qquad
q_Y^{-1}(Y^{ss}(D,T_Y)) = \b{Y}^{ss}(\b{D}, T),
$$
for any ample $D$. Consequently, the 
morphism of~\ref{rem:gelmac}~(ii)
comparing the $G_X$-quotient with the $T_Y$-quotient 
is an isomorphism.
Obviously, the family of these comparing morphisms 
is compatible with the respective GIT-systems, 
and thus defines an isomorphism of their inverse 
limits.
\end{proof}

As an immediate consequence we obtain the following rather 
special looking statements, which however give back 
the known Gelfand-MacPherson type isomorphisms of 
GIT-limits.

\begin{corollary}
Consider a $G_X$-variety $X$ and a 
$T_Y$-variety $Y$ as in~\ref{rem:gelmac}
and both projective.
If on $X$ and $Y$ every effective 
divisor is semiample, then $X$ and $Y$ have
isomorphic GIT-limits.
\end{corollary}

This applies to the case that $X$ as well as $Y$ 
are products of projective varieties having 
free divisor class group of rank one.
In particular, it applies  to the 
following setting.

\begin{corollary}
\label{cor:gelmacex}
Suppose that in the setting of~\ref{rem:gelmac},
we have $X = \PP(V_1) \times \ldots \times \PP(V_r)$ 
and $Y = \Proj((\KK[V_1] \otimes \ldots \otimes \KK[V_r])^{G^s_X})$ 
with some $G_X$-modules $V_1, \ldots, V_r$.
Then $X$ and $Y$ have isomorphic GIT-limits.
\end{corollary}

We conclude the section with a couple of examples.
The first one shows that the classical Gelfand-MacPherson 
correspondence gives rise to an isomorphism of GIT-limits;
this was observed by Kapranov~\cite{Ka}.
The second one is an analogous statement in the 
setting of complete collineations; this result 
is due to Thaddeus~\cite{Tha2}.

\begin{example}
\label{classicSL}
Consider the product $X = (\PP^{n-1})^m$, 
where $m \ge n$,
with the diagonal action of $G_X = \SL(n)$.
The total coordinate space and a Cox construction 
of $X$ are given by 
$$ 
\b{X} \ = \ (\KK^n)^m, 
\quad
\rq{X} \ = \ 
\{(v_1, \ldots, v_m)  \in \b{X}; \; 
v_i \ne 0 \text{ for } 1 \le i \le m\},
\quad
H_X \ = \ (\KK^*)^m.
$$ 
The $G_X$-action canonically lifts to 
the total coordinate space,
we have $(H_X \times G_X)^s = G_X$, and
the algebra of $G_X$-invariants is generated 
by the $(n \times n)$-minors of the matrices 
$(v_1, \ldots, v_m) \in \b{X}$.

Thus, $\b{Y} = \b{X} \quot G_X$ is 
the cone over the Grassmannian $Y := G(n,m)$,
and the basic $G_X$-invariants give 
Pl\"ucker coordinates on $\b{Y}$.
Moreover, we obtain a Cox construction 
$\rq{Y} \to Y = \rq{Y}  \quot H_Y$
with a one-dimensional torus
$H_Y \subseteq T = H_X$ and 
$\rq{Y} := \b{Y} \setminus \{\pi(0)\}$,
where $\pi \colon \b{X} \to \b{Y}$ 
is the quotient morphism.

The situation fits into Corollary~\ref{cor:gelmacex}
and thus we obtain that the actions of $G_X$
on $X$ and $T/H_Y$ on $Y$ have isomorphic
GIT-limits.
\end{example}

\begin{example}
For finite dimensional vector spaces $U,V,W$, 
consider the action of 
$G_X := \SL(U)$ on the product
\begin{eqnarray*}
X 
& := & 
\PP(\Hom(U,V)) \times \PP(\Hom(U,W)).
\end{eqnarray*}
This action lifts canonically to the total
coordinate space $\b{X} = \Hom(U,V)) \times \Hom(U,W)$,
and $\b{Y}$ is the cone over the Grassmannian 
$G(n,V \oplus W)$, where $n = \dim(U)$, acted on 
by the two dimensional torus $T = H_X$.

Take $H_Y \subseteq H_X$ such that 
$Y = \rq{Y} \quot H_Y$ is the  Grassmannian 
$G(n,V \oplus W)$. 
Then the action of $T_Y \cong \KK^*$ 
on $Y$ comes from letting act $\KK^*$ with weight 1
on $V$ and weight $-1$ on $W$.
By Corollary~\ref{cor:gelmacex}, the 
action of $G_X$ on $X$ and $T_Y$ on $Y$ have 
isomorphic GIT-limits.
\end{example}


\section{Geometry of (many) quotient spaces}
\label{sec:quotgeo}

Let a connected reductive group $G$ act 
on a normal variety $X$ 
with finitely generated Cox ring.
In this section, we show that the 
description of good $G$-sets $U \subseteq X$ 
in terms of orbit cones also opens
an approach to study the geometry 
of the quotient spaces $U \quot G$;
it turns out that in many cases the language of 
bunched rings developed in~\cite{BeHa1}
can be applied.

First, we recall the basic concepts of~\cite{BeHa1}.
Consider a factorial, 
finitely generated $\KK$-algebra $R$, 
graded by some lattice $K \cong \ZZ^{k}$
such that $R^* = \KK^*$ holds.
The latter condition enables us to 
fix a system 
$\mathfrak{F} =\{f_{1}, \dots, f_{r}\} \subset R$ 
of homogeneous pairwise non associated 
nonzero prime generators for $R$.

The {\em projected cone $(E \topto{Q} K, \gamma)$
associated\/} to the system of generators 
$\mathfrak{F} \subset R$ consists of the surjection 
$Q$ of the lattices $E := \ZZ^{r}$ and $K$ sending 
the $i$-th canonical base vector $e_{i} \in \ZZ$ 
to the degree $\deg(f_{i}) \in K$, 
and the cone $\gamma \subset  E_{\QQ}$ generated by 
$e_{1}, \dots, e_{r}$.
\begin{enumerate}    
\item
We call $\mathfrak{F} \subset R$ is {\em admissible},
if, for each facet $\gamma_{0} \preceq \gamma$,
the image $Q(\gamma_{0} \cap E)$ generates 
the lattice $K$.\item \label{fface}    
A face $\gamma_{0} \preceq \gamma$ is called an
{\em $\mathfrak{F}$-face\/} if the product
over all $f_{i}$ with $e_{i} \in \gamma_{0}$
does not belong to the ideal
$\sqrt{\bangle{f_{j}; \; e_{j} \not\in \gamma_{0}}} \subset R$.
\item \label{fbunch}
An {\em $\mathfrak{F}$-bunch\/} is  a nonempty 
collection $\Phi$ of projected $\mathfrak{F}$-faces 
with the following properties:
\begin{itemize}
\item 
a projected $\mathfrak{F}$-face $\tau$ belongs to $\Phi$ if 
and only if for each $\tau \neq \sigma \in \Phi$ we have
$\emptyset \neq \tau^{\circ} \cap \sigma^{\circ} \neq 
\sigma^{\circ}$,
\item 
for each facet $\gamma_{0} \prec \gamma$,
there is a cone $\tau \in \Phi$ such that 
$\tau^{\circ} \subseteq Q(\gamma_{0})^{\circ}$ 
holds.
\end{itemize}
\end{enumerate}
Given an $\mathfrak{F}$-bunch  $\Phi$ 
in the projected cone $(E \topto{Q} K, \gamma)$ 
associated to an admissible sytem of generators
$\mathfrak{F} \subset R$ as above, 
we call the triple $(R,\mathfrak{F},\Phi)$
a {\em bunched ring}.

Given a bunched ring $(R,\mathfrak{F},\Phi)$,
with corresponding projected cone 
$(E \topto{Q} K, \gamma)$,
consider the affine variety $Z := \Spec(R)$,
the torus $H := \Spec(\KK[K])$, and the action
$H \times Z \to Z$ given by the $K$-grading of
$R$.
Then the projected $\mathfrak{F}$-faces
are precisely the orbit cones of the 
$H$-action on $Z$, and there is a canonical injection
\begin{eqnarray*}
\{\mathfrak{F}  \text{-bunches} \}
& \rightarrow &
\{\text{2-maximal collections in } \Omega(Z)\}
\\
\Phi
& \mapsto &
\Psi(\Phi) \ := \ 
\left\{
\omega(z); \; 
z \in Z, \; 
\tau^\circ \subset \omega(z)^\circ 
\text{ for some } \tau \in \Phi
\right\}.
\end{eqnarray*}
Using this observation, we may associate
to the bunched ring 
$(R,\mathfrak{F},\Phi)$
a variety by setting
\begin{eqnarray*}
X(R,\mathfrak{F},\Phi)
& = & 
U(\Psi(\Phi)) \quot H.
\end{eqnarray*}

The main object of~\cite{BeHa1} is to
read off geometric properties of this 
variety from its defining data.
Let us briefly provide the necessary
notions.
Call an $\mathfrak{F}$-face 
$\gamma_{0} \preceq \gamma$ {\em relevant\/}
if $Q(\gamma_{0})^{\circ} \supset \tau^{\circ}$ 
holds for some $\tau \in \Phi$,
and denote by $\rlv(\Phi)$ the collection
of relevant $\mathfrak{F}$-faces.
The {\em covering collection} of $\Phi$ 
is the collection $\cov(\Phi) \subset \rlv(\Phi)$
of set-theoretically minimal members of 
$\rlv(\Phi)$.

\begin{theorem}
\label{bunchgeom}
Let $(R,\mathfrak{F},\Phi)$ be a bunched ring
with corresponding projected cone 
$(E \topto{Q} K, \gamma)$,
and let $X := X(R,\mathfrak{F},\Phi)$ 
be the associated variety.
\begin{enumerate}
\item 
The variety 
$X$ is locally factorial if and only if
$Q(\gamma_0 \cap E)$ generates the lattice 
$K$ for every $\gamma_0 \in \rlv(\Phi)$.
\item 
The variety 
$X$ is $\QQ$-factorial if and only if 
any cone of $\Phi$ is of full dimension 
in $K_\QQ$. 
\item 
The dimension of $X$ is $\dim(R)-\dim(K_\QQ)$,
its divisor class group  is 
$\Cl(X) \cong K$, and the Picard group of 
$X$ sits in $\Cl(X)$ as 
\begin{eqnarray*}
\Pic(X) 
& = & 
\bigcap_{\gamma_{0} \in \cov(\Phi)} Q(\lin(\gamma_{0}) \cap E).
\end{eqnarray*}
\item 
The effective cone, the moving cone,
and the cones of semiample and ample 
divisor classes of X are given by
$$
\begin{array}{ccc}
\displaystyle \Eff(X) \; = \; Q(\gamma),
& \qquad &
\displaystyle \Mov(X) \; = \; 
\bigcap_{\gamma_0 \text{ facet of } \gamma} Q(\gamma_0),
\\[3ex]
\displaystyle \SAmple(X) \; = \; \bigcap_{\tau \in \Phi} \tau,
& \qquad &
\displaystyle \Ample(X) \; = \; \bigcap_{\tau \in \Phi} \tau^{\circ}.
\end{array}
$$
\item
Suppose that for $d := r - \dim(R)$,
we have  $K$-homogeneous generators  $g_{1}, \ldots, g_{d}$ 
for the relations between the members
$f_{1}, \ldots, f_{r}$ of $\mathfrak{F}$. 
Then the canonical divisor class of $X$ is
$$ 
\sum \deg(g_{j}) - \sum \deg(f_{i}).
$$
\end{enumerate}
\end{theorem}

Now consider a factorial affine algebra
$R$ graded by a lattice $K$, let 
$\mathfrak{F} \subset R$ be an admissible system of 
generators,
and let $(E \topto{Q} K, \gamma)$ be the associated
projected cone. 
Consider the affine variety $Z := \Spec(R)$
and the action of the torus $H := \Spec(\KK[K])$ 
on it given by the $K$-grading of $R$. 
Define a subset of the weight cone $\omega(Z)$ by
$$
\omega(R)^{\odot}
\ := \ 
\omega(Z)^{\odot}
\ := \
\bigcap_{\gamma_{0} \preceq \gamma} Q(\gamma_0)^0,  
$$
where $\gamma_0$ runs through all facets of $\gamma$. 
With any $\chi\in\omega(Z)$ associate the collection
$\Psi(\chi)$ of projected cones $\tau$ with $\chi\in\tau^0$ and
the collection $\Phi(\chi)$ of set-theoretically minimal elements
of $\Psi(\chi)$.

\begin{proposition}
\label{prop:deepint}
Suppose that
for each facet $\gamma_{0} \preceq \gamma$,
the image $Q(\gamma_{0} \cap E)$ generates 
the lattice $K$ and 
that $R_0=\KK$ holds. 
Then for any $\chi\in\omega(Z)^{\odot}$ the triple
$(R,\mathfrak{F},\Phi(\chi))$ is the bunched ring representing
the quotient space $Z^{ss}(\chi)\quot H$.
\end{proposition}

\begin{proof}
The condition $R_0=\KK$ implies that $Z^{ss}(\chi)\quot H$ 
is projective 
and thus the collection $\Psi(\chi)$ is 2-maximal. 
This shows that $\Phi(\chi)$
satisfies the first condition in the definition of 
an $\mathfrak{F}$-bunch. Note also that for any
facet $\gamma_0\preceq\gamma$ the cone $Q(\gamma_0)$ 
is an orbit cone. 
Indeed, since $\mathfrak{F}$ is a system of 
pairwise non associated nonzero prime generators of $R$,
the condition $f_i(z)=0$ does 
not imply $f_j(z)=0$ for any $j\ne i$. 
Since $Q(\gamma_0)\in\Psi(\chi)$, 
there is an element $\tau\in\Phi(\chi)$ with
$\tau^0\subseteq Q(\gamma_0)^0$. 
We have checked that $\Phi(\chi)$ is an 
$\mathfrak{F}$-bunch.
The other statements follow from the definition 
of variety associated with a bunched ring.
\end{proof}

\begin{corollary}
\label{amplechamber}
For any $\chi\in\omega(Z)^{\odot}$ the cone $\omega(Z)$ is the
cone of effective divisors and the GIT-cone $\lambda(\chi)$ is the
cone of semiample divisors of the variety $Z^{ss}(\chi)\quot H$.
\end{corollary}

Now, let a connected reductive group $G$ act
on a normal projective variety $X$.
Suppose that there is a $G$-equivariant Cox 
construction $\rq{X} \to X = \rq{X} \quot H_X$
with total coordinate space $\b{X}$.
Then the invariant Cox ring 
$R := \mathcal{R}(X)^{G^s}$ comes with a grading
by the character group $K$ of $H := H_X \times G/G^s$ 
corresponding to the induced $H$-action
on $\b{X} \quot G^s$,
and Proposition~\ref{prop:deepint} 
gives the following.

\begin{corollary}
Let $\mathfrak{F}$ be an admissible system
of generators for $R = \mathcal{R}(X)^{G^s}$.
Then for any $G$-linearized ample divisor $D$ on $X$ 
defining an element $\chi \in \omega(R)^{\odot}$,
the associated quotient space $X^{ss} \quot G$ 
arises from the bunched ring 
$(R,\mathfrak{F},\Phi(\chi))$.
\end{corollary}

In a first example, we consider once 
more the diagonal action of the special 
linear group $\SL(n)$ on a product of 
projective spaces $\PP^{n-1}$.
It has quite a big variation of 
GIT-quotients, but there is one 
``canonical'' candidate, namely the
(unique) set of semistable points, 
which is invariant under permuting 
the factors  $\PP^{n-1}$.
The case $n =2$, is studied by several
authors, see~\cite{BBSo} for a 
uniqueness result and~\cite{Po} for
an approach to the geometry.
We will see, also for higher $n$, that
the quotient fits into the setting
of bunched rings.

\begin{example}
Consider the diagonal action of
$G=\SL(n)$ on $X = (\PP^{n-1})^m$, 
where $m \ge n+2$.
We identify $\ZZ^m \to \Cl(X)$ by 
sending $e_i$ to the class of $D_i$ 
of the pullback of a hyperplane in 
the $i$-th factor $\PP^{n-1}$.
Then the action of  $H = (\KK^*)^m$
on the total coordinate space 
$\b{X} = (\KK^n)^m$ is 
componentwise scalar multiplication.

The quotient $\b{Y} = \b{X} \quot G$
is the cone over the Grassmannian 
$G(n,m)$,
and $R := \KK[\b{Y}]$ is generated 
by the $(n \times n)$-minors of the matrices 
$(v_1, \ldots, v_m) \in \b{X}$.
The weights of these generators with respect
to the $H$-action 
are the $w \in \ZZ^m$ having exactly $n$
entries 1 and the others 0.
Consequently, the weight cone is
$$
\omega(\b{X})
\ = \ 
\omega(\b{Y})
= 
\left\{
(s_1, \ldots, s_m) \in \QQ^m_{\ge 0}; \;
s_1 + \ldots + s_m \ge ns_j,  \, 1 \le j \le m
 \right\}.
$$
The associated GIT-fan is known, see~\cite[3.3.24]{DoHu}.
It is most conveniently described by giving the
walls; these are the intersections of 
$\omega(\b{X})$ with the hyperplanes
\begin{eqnarray*}
\mathcal{H}_{k,J}
& := & 
\left\{
(s_1,\dots,s_m); \; 
(n-k)\sum_{j \in J} s_j\ = \ k\sum_{l \notin J} s_l
\right\},
\\
& & 
1\le k\le\frac{n}{2}, \quad 
J\subset \{1,\dots,m\}, 
\quad
k<|J|<m+k-n.
\end{eqnarray*}
The $(n \times n)$-minors mentioned before
form a system of pairwise nonassociated
prime generators for $R$, 
and one easily checks 
that we are in the situation 
of Proposition~\ref{prop:deepint}.
Thus, in order to figure out the GIT-quotients, for 
which we get a describing bunched ring for 
free, we need the cone $\omega(\b{X})^{\odot}$;
it is given as 
$$
\left\{
(s_1, \ldots, s_m) \in \QQ^m_{\ge 0}; \
\sum_{j\in J} s_j < (n-1)\sum_{l\notin J} s_l, \ 
 J \subset \{1,\dots,m\}, \ 
|J|=n
\right\}.
$$ 

There is a unique set $U \subseteq X$ 
of semistable points, which is invariant
under permutation of the factors $\PP^{n-1}$ 
of $X$; it is defined by the divisor 
$D:=D_1+\dots+D_m$.
In our picture, the class of  $D$ is the 
point $\chi := (1,\dots,1)$.
The corresponding GIT-cone $\lambda(\chi)$ 
is defined by inequalities
\begin{eqnarray*}
(n-k)\sum_{j\in J} s_j \le k\sum_{l\notin J} s_l, 
&  & 
1\le k\le\frac{n}{2}, \
J\subset \{1,\dots,m\}, \
k \le |J| \le \frac{km}{n},
\\
k\sum_{l\notin J} s_l \le (n-k)\sum_{j\in J} s_j, 
& & 
1\le k\le\frac{n}{2}, \
J\subset \{1,\dots,m\}, \
\frac{km}{n} \le |J| \le m+k-n.
\end{eqnarray*}
Note that here only the cases  
$km/n-1 < |J| \le km/n+1$ are essential. 
In particular, we see that 
$\chi$ belongs to $\omega(\b{X})^{\odot}$ 
provided $m\ge 5$ and $n=2$, or $m\ge n+2$ and $n\ge 3$. 
From Theorem~\ref{bunchgeom} we infer that $\lambda(\chi)$
is the semiample cone of the quotient space $U \quot G$.
Let us have a closer look at it.

\smallskip

\noindent {\em Case 1. } We have $\text{GCD}(m,n) = d>1$. 
Then $\chi$ lies in all walls corresponding to 
$k = n/d$ and $|J|= m/d$.
It follows that the GIT-cone $\lambda(\chi)$ 
is a ray. In particular, the $U \quot G$ comes
with non-$\QQ$-factorial singularities, and its
Picard group is of rank one.
 
\smallskip

\noindent
{\em Case 2. } The numbers $m$ and $n$ are coprime. 
Then $\chi$ is not contained in any wall, 
so $\lambda(\chi)$ has full dimension.
In the cases $n=2,3$ we obtain the following 
describing inequalities ($S_m$ stands for the
symmetric group). 
\begin{itemize}
\item 
For $n=2$ and $m=2r+1$, the cone  $\lambda(\chi)$
is given by the inequalities
\begin{eqnarray*}
s_{\pi(1)} + \ldots + s_{\pi(r)} 
& \le &   s_{\pi(r+1)} + \ldots + s_{\pi(m)}, 
\qquad
\pi \in S_m.
\end{eqnarray*}
\item
For $n=3$ and $m=3r+s$, $s=1,2$,  the cone  $\lambda(\chi)$
is given by 
\begin{eqnarray*}
2(s_{\pi(1)} + \ldots + s_{\pi(r)}) 
& \le & 
s_{\pi(r+1)} + \ldots + s_{\pi(m)}, 
\qquad
\pi \in S_m,
\\
2(s_{\pi(1)} + \ldots + s_{\pi(r+1)}) 
& \ge & 
s_{\pi(r+2)} + \ldots + s_{\pi(m)}, 
\qquad
\pi \in S_m.
\end{eqnarray*}
\end{itemize}

The number of extremal rays of the semiample cone is
certainly an invariant for any variety;
here we obtain, for example, that for 
$n=3$ and $m=5$ the quotient $U \quot G$ 
is a smooth projective surface having a semiample 
cone with 10 extremal rays.
\end{example}

Also for the symplectic group action on a product
of projective spaces studied in Theorem~\ref{tsp}, 
we have a unique set of semistable points being
invariant under permuting the factors.
Here comes more information.

\begin{example}
Consider the diagonal action of 
the symplectic group $\Sp(2n)$ on $(\PP^{2n-1})^m$
as in Theorem~\ref{tsp}.
Then the cone $\omega(\b{Y})^{\odot} \subseteq \QQ^m$ 
is given by the additional inequalities
\begin{eqnarray*}
s_{\pi(1)} + s_{\pi(2)} 
& \le &   s_{\pi(3)} + \ldots + s_{\pi(m)}, 
\qquad
\pi \in S_m.
\end{eqnarray*}
In particular, the point $\chi=(1,\dots,1)$ 
is contained in $\omega(Z)^{\odot}$
for $m\ge 5$. 
Moreover, this point belongs to 
the walls $K_{J_1,J_2}$ with $|J_1|=|J_2|$. 
Taking $|J_1|=|J_2|=1$, 
one gets that the GIT-cone 
$\lambda(\chi)$ is one-dimensional.
So the resulting quotient space has 
non-$\QQ$-factorial singularities
and its Picard group is of rank one.
\end{example}

The next example belongs to a large class,
arising from reducible $G$-representations  
whose algebra of invariants has a single 
relation. 
A complete classification for $G=\SL(n)$ 
is given in~\cite{Shm}. 
There also the weights of generators and 
of the relation are listed.
However, the relation itself is sometimes 
not easy to write down explicitly. 
This may cause difficulties in determining 
orbit cones.
The following observation helps.

\begin{lemma}
\label{lem2112}
Let a torus $T$ act diagonally on $\KK^r$ 
with weights $(w_1, \ldots, w_r)$, and 
consider a $T$-invariant hypersurface 
$Z := V(\KK^r;f)$ with a polynomial
$f$ of the form
$$ 
f \ = \ T_r^k + g, \text{ where } g \in \KK[T_1, \ldots, T_{r-1}].
$$
Then the orbit cones of the $T$-action 
on $Z$ are precisely the cones $\cone(w_j; \ j \in J)$ 
with $J \subseteq \{1, \ldots, r-1\}$,
and the GIT-fan of the $T$-action on 
$Z$ is $\Sigma(w_1, \ldots, w_{r-1})$.
\end{lemma}

\begin{proof}
First we show that any orbit cone 
$\omega(z)$, $z \in Z$, is of the form
$\cone(w_j; \ j \in J)$ 
with a subset $J \subseteq \{1, \ldots, r-1\}$.
If $z_r = 0$ holds, then $\omega(z)$ is 
necessarily generated by weights from
$\{w_1, \ldots, w_{r-1}\}$.
For $z_r \ne 0$, we have
$g(z_1, \ldots, z_{r-1}) \ne 0$. 
Thus, some monomial 
$g_0 = T_{i_1}^{\nu_{i_1}} \cdots T_{i_l}^{\nu_{i_l}}$ 
occuring in $g$ satisfies $g_0(z) \ne 0$.
This implies
$$
k w_r 
\ = \ 
\deg(g)  
\ \in \ 
\cone(w_{i_1}, \ldots, w_{i_l})
\ \subseteq \
\omega(z).
$$
Consequently, we see that $w_r$ is not needed 
as a generator of the orbit cone $\omega(z)$, 
and we are done.

Conversely, let $\omega = \cone(w_j; \ j \in J)$ 
with a subset $J \subseteq \{1, \ldots, r-1\}$
be given. Then we need a point $z \in Z$
with $\omega(z) = \omega$. 
For $1 \le j \le r-1$, we set
\begin{eqnarray*}
z_j 
& := &
\begin{cases}
1 & \text{for } j \in J,
\\
0 & \text{for } j \not\in J.
\end{cases}
\end{eqnarray*}
Next take $z_r \in \KK$ with $z_r^k = -g(z_1, \ldots, z_{r-1})$
Then $z := (z_1, \ldots, z_r)$ belongs to $Z$ 
because of $f(z) = 0$, and we directly see
$\omega \subseteq \omega(z)$.

For $z_r = 0$, also the inclusion 
$\omega(z) \subseteq \omega$ is obvious.
If $z_r \ne 0$ holds, then we have 
$g(z_1, \ldots, z_{r-1})  \ne  0$,
and a similar reasoning as before gives
$w_r \in \omega$.
Moreover, if $h$ is any semiinvariant with $h(z) \ne 0$,
then some monomial
$h_0 = T_{i_1}^{\nu_{i_1}} \cdots T_{i_l}^{\nu_{i_l}}$ 
satisfies $h_0(z) \ne 0$.
The latter condition implies 
$i_1, \ldots, i_l \in J \cup \{r\}$,
which in turn gives $\deg(h) \in \omega$.
\end{proof}

\begin{example}
Let $G := \SL(4)$ and consider the irreducible 
representations $G \to \GL(V)$ and 
$G \to \GL(W)$ with the respective 
highest weights $\omega_2$ and $\omega_2^2$, 
where  $\omega_2$ denotes the second fundamental 
weight, i.e., we have $V := \bigwedge^2 \KK^4$
and $W$ is a 20-dimensional subspace of $S^2V$.
Then we have an induced $G$-action on 
\begin{eqnarray*}
X 
& := & 
\PP(V) \times \PP(W).
\end{eqnarray*} 
By construction,
$\b{X} = V \times W$ is the equivariant total 
coordinate space.
The Neron-Severi
torus $H = (\KK^*)^2$ acts by componentwise scalar multiplication 
on $\b{X}$.
According to~\cite[Table~4]{Shm}, the algebra of invariants is
of the form
$$ 
R \ :=\ 
\KK[\b{X}]^G
\ = \ 
\KK[T_1,\dots,T_{12}] / \bangle{f},
\qquad
f = T_{12}^2 + g(T_1, \ldots, T_{11})
$$
The classes $f_i \in R$ of the $T_i$ form 
an admissible system 
$\mathfrak{F} = (f_1, \ldots, f_{12})$
of pairwise nonassociated prime homogeneous
generators.
Their degrees $w_1 \ldots, w_{12}$ were 
also calculated in~\cite{Shm}; in $\ZZ^2$ 
they can be given as
$$
w_{i} =  (2+i,0), \  i=0,\ldots,4,
\qquad
w_{6+j} = (j,1),  \ j=0,\ldots,5,
\qquad
w_{12} = (15,3).
$$

According to Lemma~\ref{lem2112}, the GIT-fan 
of the $(H \times G)$-action on $\b{X}$ is 
$\Sigma(w_1, \ldots, w_{11})$.
The situation is sketched in the
following figure, where 
the  bullets stand for the weights, 
and the shadowed area indicates the 
cone $\omega(\b{X})^\odot$. 
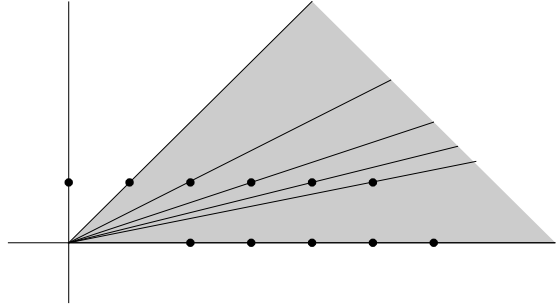
\begin{center}
\input{shmelkinex2.pstex_t}
\end{center}
Thus, we have five $\QQ$-factorial 
and four non-$\QQ$-factorial quotient
spaces, to which Theorem~\ref{bunchgeom}
directly can be applied.
Note that one of them is $\QQ$-Fano (but not Fano), 
and each of them comes with singularities.
\end{example}


\end{document}

%% file: butterfly1.pstex_t
\begin{picture}(0,0)%
\includegraphics{butterfly1.pstex}%
\end{picture}%
\setlength{\unitlength}{1243sp}%
\begingroup\makeatletter\ifx\SetFigFont\undefined%
\gdef\SetFigFont#1#2#3#4#5{%
  \reset@font\fontsize{#1}{#2pt}%
  \fontfamily{#3}\fontseries{#4}\fontshape{#5}%
  \selectfont}%
\fi\endgroup%
\begin{picture}(3666,3666)(2668,-3694)
\end{picture}%

%% file: butterfly3.pstex_t
\begin{picture}(0,0)%
\includegraphics{butterfly3.pstex}%
\end{picture}%
\setlength{\unitlength}{1243sp}%
\begingroup\makeatletter\ifx\SetFigFont\undefined%
\gdef\SetFigFont#1#2#3#4#5{%
  \reset@font\fontsize{#1}{#2pt}%
  \fontfamily{#3}\fontseries{#4}\fontshape{#5}%
  \selectfont}%
\fi\endgroup%
\begin{picture}(4084,4084)(2475,-3662)
\put(3601,-1411){\makebox(0,0)[lb]{\smash{{\SetFigFont{7}{8.4}{\familydefault}{\mddefault}{\updefault}{\color[rgb]{0,0,0}$\omega_3^+$}%
}}}}
\put(5401,-736){\makebox(0,0)[lb]{\smash{{\SetFigFont{7}{8.4}{\familydefault}{\mddefault}{\updefault}{\color[rgb]{0,0,0}$\omega_2^+$}%
}}}}
\put(4276,-2761){\makebox(0,0)[lb]{\smash{{\SetFigFont{7}{8.4}{\familydefault}{\mddefault}{\updefault}{\color[rgb]{0,0,0}$\omega_1^+$}%
}}}}
\end{picture}%

%% file: butterfly2.pstex_t
\begin{picture}(0,0)%
\includegraphics{butterfly2.pstex}%
\end{picture}%
\setlength{\unitlength}{1243sp}%
\begingroup\makeatletter\ifx\SetFigFont\undefined%
\gdef\SetFigFont#1#2#3#4#5{%
  \reset@font\fontsize{#1}{#2pt}%
  \fontfamily{#3}\fontseries{#4}\fontshape{#5}%
  \selectfont}%
\fi\endgroup%
\begin{picture}(4084,4084)(2668,-3919)
\put(5626,-1861){\makebox(0,0)[lb]{\smash{{\SetFigFont{7}{8.4}{\familydefault}{\mddefault}{\updefault}{\color[rgb]{0,0,0}$\omega_3^-$}%
}}}}
\put(4726,-736){\makebox(0,0)[lb]{\smash{{\SetFigFont{7}{8.4}{\familydefault}{\mddefault}{\updefault}{\color[rgb]{0,0,0}$\omega_1^-$}%
}}}}
\put(3376,-2761){\makebox(0,0)[lb]{\smash{{\SetFigFont{7}{8.4}{\familydefault}{\mddefault}{\updefault}{\color[rgb]{0,0,0}$\omega_2^-$}%
}}}}
\end{picture}%

%% file: schwarzex.pstex_t
\begin{picture}(0,0)%
\includegraphics{schwarzex.pstex}%
\end{picture}%
\setlength{\unitlength}{1657sp}%
\begingroup\makeatletter\ifx\SetFigFont\undefined%
\gdef\SetFigFont#1#2#3#4#5{%
  \reset@font\fontsize{#1}{#2pt}%
  \fontfamily{#3}\fontseries{#4}\fontshape{#5}%
  \selectfont}%
\fi\endgroup%
\begin{picture}(4260,5023)(1876,-5201)
\put(1891,-3661){\makebox(0,0)[lb]{\smash{{\SetFigFont{8}{9.6}{\familydefault}{\mddefault}{\updefault}{\color[rgb]{0,0,0}$w_5$}%
}}}}
\put(2251,-4111){\makebox(0,0)[lb]{\smash{{\SetFigFont{8}{9.6}{\familydefault}{\mddefault}{\updefault}{\color[rgb]{0,0,0}$w_4$}%
}}}}
\put(4231,-2401){\makebox(0,0)[lb]{\smash{{\SetFigFont{8}{9.6}{\familydefault}{\mddefault}{\updefault}{\color[rgb]{0,0,0}$w_3$}%
}}}}
\put(6121,-5101){\makebox(0,0)[lb]{\smash{{\SetFigFont{8}{9.6}{\familydefault}{\mddefault}{\updefault}{\color[rgb]{0,0,0}$w_2$}%
}}}}
\put(1891,-421){\makebox(0,0)[lb]{\smash{{\SetFigFont{8}{9.6}{\familydefault}{\mddefault}{\updefault}{\color[rgb]{0,0,0}$w_1$}%
}}}}
\put(1891,-1861){\makebox(0,0)[lb]{\smash{{\SetFigFont{8}{9.6}{\familydefault}{\mddefault}{\updefault}{\color[rgb]{0,0,0}$w_6$}%
}}}}
\end{picture}%

%% file: shmelkinex2.pstex_t
\begin{picture}(0,0)%
\includegraphics{shmelkinex2.pstex}%
\end{picture}%
\setlength{\unitlength}{1657sp}%
\begingroup\makeatletter\ifx\SetFigFontNFSS\undefined%
\gdef\SetFigFontNFSS#1#2#3#4#5{%
  \reset@font\fontsize{#1}{#2pt}%
  \fontfamily{#3}\fontseries{#4}\fontshape{#5}%
  \selectfont}%
\fi\endgroup%
\begin{picture}(8124,4524)(1339,-4123)
\end{picture}%

%% file: gquots.bbl
\begin{thebibliography}{}%
%
\bibitem{AG} O.M.~Adamovich, E.O.~Golovina: 
Simple linear Lie groups having a free algebra of invariants.
In: Voprosy Teorii Grupp i Gomologicheskoj Algebry 2, 3--41 (1979). 
English Transl.: Sel. Math. Sov. 3:2, 183-220 (1984) 
%
\bibitem{ArHa} I.V.~Arzhantsev, J.~Hausen: 
On embeddings of homogeneous spaces with small boundary. 
J.~Algebra~304, No.~2, 950-988 (2006),
{\tt math.AG/0507557}
%
\bibitem{BeHa1} F.~Berchtold, J.~Hausen: 
Cox rings and combinatorics.
Transactions of the AMS~359, No.~3, 
1205-1252 (2007),
{\tt math.AG/0311105}
%
\bibitem{BeHa2} F.~Berchtold, J.~Hausen: 
GIT-equivalence beyond the ample cone, 
Michigan Math. J.~54, No.~3, 483-516 (2006), 
{\tt math.AG/0503107}
%
\bibitem{BB2} A.~Bia\l ynicki-Birula: Algebraic
  Quotients. In: R.V.~Gamkrelidze, V.L.~Popov (Eds.), Encyclopedia of
  Mathematical Sciences, Vol. 131., 1--82 (2002) 
%
\bibitem{BBSo}  A.~Bia\l ynicki-Birula, A.J.~Sommese: 
Quotients by $\CC^*$ and $\SL(2,\CC)$ actions.
Trans. Amer. Math. Soc. 279, 
773--800 (1983)
%
\bibitem{BBSw4}  A.~Bia\l ynicki-Birula, J.~\'Swi\c{e}cicka: Three
  theorems on existence of good quotients. Math. Ann., Vol.~307,
  143--149 (1997)
%
\bibitem{BBSw5}  A.~Bia\l ynicki-Birula, J.~\'Swi\c{e}cicka: 
A recipe for finding open subsets of vector spaces with a good quotient.  
Colloq. Math.  77,  97--114  (1998)
%
\bibitem{Co} D.~Cox: The homogeneous coordinate ring of a toric
  variety. J. Alg. Geom. 4, 17--50  (1995) 
%
\bibitem{DoHu} I.V.~Dolgachev, Y.~Hu: Variation of geometric 
invariant theory quotients. (With an appendix: ``An example of a 
thick wall'' by N.~Ressayre).
Publ. Math., Inst. Hautes Etud. Sci. 87, 5--56 (1998)
%
\bibitem{GM} I.M.~Gelfand, R.W.~MacPherson: 
Geometry in Grassmanians and a generalization of the dilogarithm.
Adv. in Math. 44, 279--312 (1982)
%
\bibitem{EGA}
A.~Grothendieck:
\'{E}l\'{e}ments de g\'{e}om\'{e}trie alg\'{e}brique. 
Publ. Math., Inst. Hautes \'{E}tud. Sci. 4, 
1-228 (1960); 8, 1-222 (1961); 11, 349-511 (1962)
%
\bibitem{Ha} J.~Hausen: Geometric Invariant Theory based on Weil divisors.
     Compositio Math. 140,  1518-1536 (2004)
%
\bibitem{Har}
R.~Hartshorne: Algebraic Geometry, 
GTM~52, Springer Verlag (1977)
%
\bibitem{Kn}
F.~Knop:
\"Uber Hilberts vierzehntes Problem 
f\"ur Variet\"aten mit Kompliziertheit eins. 
Math. Z. 213, 33--36 (1993)
%
\bibitem{Lu} D.~Luna: 
Sur les orbites ferm\'ees des groupes alg\'ebriques reductif.
Invent. Math. 16, 1--5 (1972)
%
\bibitem{Mu} D.~Mumford, J.~Fogarty, F.~Kirwan: Geometric 
  Invariant Theory. 3rd enl. ed.. Ergebnisse der 
  Mathematik und ihrer Grenzgebiete. Berlin: Springer-Verlag.
  (1993)
%
\bibitem{Po} M.~Polito: 
$SL(2,\CC)$-quotients de $(\PP^1)^n$. C. R. Acad. Sci. Paris,
t.~321, S\'erie I, 1577--1582 (1995)
%
\bibitem{PV2} V.L.~Popov, E.B.~Vinberg: 
Invariant theory; 
Itogi Nauki i Tekhniki, Sovrem. Probl. Mat. Fund. Naprav. 
Vol.~55, VINITI, Moscow, 137--314 (1989). 
English Transl.: Algebraic Geometry IV, Encyclopedia 
of Math. Sciences, vol.~55, Springer-Verlag, Berlin (1994)
%
\bibitem{Ka} M.M.~Kapranov: 
Chow quotients of Grassmannians. I.
Advances in Soviet Math. 16, Part 2, 29--110 (1993)
%
\bibitem{Ki} A.~King: 
Moduli of representations of finite-dimensional algebras.  
Quart. J. Math. Oxford Ser. (2)  45,  
515--530  (1994)
%
\bibitem{Re} N.~Ressayre: 
The GIT-equivalence for $G$-line bundles. 
Geom. Dedicata 81, No. 1--3, 295--324 (2000)
%
\bibitem{Sch} G.W.~Schwarz: 
Representations of simple Lie groups with regular rings of invariants.
Invent. Math. 49, 167--191 (1978)
%
\bibitem{Se}
C.S.~Seshadri:
Quotient spaces modulo reductive algebraic groups.  
Ann. of Math.~(2)  95, 511--556  (1972)
%
\bibitem{Shm} D.A.~Shmel'kin: 
On representations of $SL(n)$ with algebras of invariants 
being complete intersections.
J. Lie Theory 11, 207--229 (2001)
%
\bibitem{Tha} M.~Thaddeus: Geometric invariant theory and flips.  
J.~Amer.~Math.~Soc.~9, 691--723 (1996)
%
\bibitem{Tha2} M.~Thaddeus: Complete collineations revisited.
Math. Ann.~315, 469--495 (1996)
%
\bibitem{Wl} J.~W\l odarczyk: 
Embeddings in toric varieties and prevarieties.  
J. Algebraic Geom.~2, 705--726 (1993)
\end{thebibliography}
